\documentclass[oneside,12pt]{amsart}
\usepackage{amscd,amsmath,latexsym, amssymb}
\usepackage{xcolor}
\usepackage{tikz}
\usepackage{enumitem}
\usepackage{dsfont}
\usepackage{amsthm}
\usepackage{mathtools}
\usetikzlibrary{matrix}
\usepackage{verbatim}
\usepackage{amsfonts}
\newtheorem{thm}[equation]{Theorem}
\theoremstyle{plain}
\newtheorem*{theorem*}{Theorem}
\numberwithin{equation}{section}
\newtheorem{cor}[equation]{Corollary}
\newtheorem{exm}[equation]{Example}
\newtheorem{rmk}[equation]{Remark}

\newtheorem{lem}[equation]{Lemma}

\newtheorem{defin}[equation]{Definition}

\newtheorem{prop}[equation]{Proposition}

\newtheorem{assum}[equation]{Assumption}

\newtheorem{alg}[equation]{Algorithm}

\begin{document}
%\raggedbottom \voffset=-.7truein \hoffset=0truein \vsize=8truein
%\hsize=6truein \textheight=8truein \textwidth=6truein
\baselineskip=18truept
\def\mapright#1{\ \smash{\mathop{\longrightarrow}\limits^{#1}}\ }
\def\mapleft#1{\ \smash{\mathop{\longleftarrow}\limits^{#1}}\ }
\def\mapup#1{\Big\uparrow\rlap{$\vcenter {\hbox {$#1$}}$}}
\def\mapdown#1{\Big\downarrow\rlap{$\vcenter {\hbox {$\ssize{#1}$}}$}}
\def\mapne#1{\nearrow\rlap{$\vcenter {\hbox {$#1$}}$}}
\def\mapse#1{\searrow\rlap{$\vcenter {\hbox {$\ssize{#1}$}}$}}
\def\mapr#1{\smash{\mathop{\rightarrow}\limits^{#1}}}
\def\lb{[}
\def\ss{\smallskip}
\def\at{{\widetilde\alpha}}
\def\sm{\wedge}
\def\la{\langle}
\def\ra{\rangle}
\def\on{\operatorname}
\def\qed{\quad\rule{8pt}{8pt}\bigskip}
\def\ssize{\scriptstyle}
\def\ar#1{\stackrel {#1}{\rightarrow}}
\def\br{\bold R}
\def\bc{\bold C}
\def\si{\sigma}
\def\zp{\bold Z_p}
\def\da{\downarrow}
\def\xbar{{\overline x}}
\def\ebar{{\overline e}}
\def\zbar{ {\overline z}}
\def\ni{\noindent}
\def\coef{\on{coef}}
\def\den{\on{den}}
\def\ot{\otimes}
\def\ms{\medskip}
\def\cl{ \mathcal{L }}
\def\cf{\mathcal{F}}
\def\ct{\mathcal{T}}
\def\ctl{\mathcal{T}_L}
\def\ch{\mathcal{H}}
\def\d{\displaystyle\mathop}
\def\wcm{ \widetilde{\cm}}
% Rick's additions
\newtheorem{expls}[equation]{Examples}
\newtheorem{exam}[equation]{Example}
\newtheorem{exams}[equation]{Examples}
\def\ia{{\item[\rm (a)]\ }}
\def\ib{{\item[\rm (b)]\ }}
\def\ic{{\item[\rm (c)]\ }}
\def\id{{\item[\rm (d)]\ }}
\def\ad{{\rm ad}}
\def\im{{\rm im}}
\def\hl{{ \hat{\ell}}}
\def\hn{{ \hat{ \mathcal N}}}
\def\ead{{\rm expad}}
\def\nz{\mathbb{Z}}
\def\nk{\mathbb{K}}
\def\nn{\mathbb{N}}
\def\nzp{{\mathbb{Z}^+}}
\def\nr{\mathbb{R}}
\def\nc{\mathbb{C}}
\def\cg{ \mathcal{G }}
\def\cn{ \mathcal{N }}
\def\cw{ \mathcal{W }}
\def\cm{ \mathcal{M }}
\def\cl{ \mathcal{L }}
\def\chapter{\mathcal{H}}
\def\lhn{ \hat{n}}
\def\pf{\noindent {\bf Proof :\ }}
\def\qed{\hfill {\bf q.e.d.}}
\def\qedalone{\hspace*{1in}\hfill ${\rm {\bf q.}{\bf e.}{\bf d.}}$}
\def\glnc{\mathfrak{gl}(n,\nc)}
\def\glnr{\mathfrak{gl}(n,R)}
\def\vsmfv{\vspace{-.05in}}
\def\vspfv{\vspace{.05in}}
\def\vspsv{\vspace{.07in}}
\def\pr{^{\,\prime}}
\def\mh{ \widehat{m}}
\def\per{.$\!$\, }
\def\fg{\mathfrak{g}}
\def\fl{\mathfrak{L}}
\def\sl{\mathfrak{sl}}
\def\gl{\mathfrak{gl}}
\def\ibul{\item[$\bullet$]}
\def\dpr{^{\prime \prime}}
\def\without{w.l.o.g\per }
\def\pgs{pp\per }
\def\skp#1{\vskip#1cm\relax}
\def\nd{\noindent}
\def\ha{H_{\ast}}
\def\x{\mathbf{x}}
\def\p{\mathbf{p}}
\def\q{\mathbf{q}}
\def\ts{\textstyle{}}
\def\mci{\mathcal{I}}
\def\mcl{\mathcal{L}}
\def\mbk{\mathbf{k}}
\def\mcm{\mathcal{M}}
\newcommand{\ds}{\displaystyle{}}
\renewcommand\qed{ \hspace{5.4 in} $\square$}

\setlength{\oddsidemargin}{.25in}
\setlength{\textwidth}{6in}

\setlength{\voffset}{.5in}

\title[Spectral Sequence] { A Spectral Sequence for a Graded Linear Map}

\author[L.~Bates, M.~Bendersky,R.C.~Churchill]{L.M.~ Bates, M.~Bendersky, R.C.  Churchill}

\address{Bates: Dept. of Math., University of Calgary.}

\address{Bendersky and Churchill:  Dept.  of Math., Hunter College/Graduate Center, CUNY.
}
\email{}



\begin{abstract}
	We apply the method of spectral sequences to study classical problems in analysis.  We illustrate the method by finding polynomial vector fields that commute with a given polynomial vector field and finding integrals of polynomial Hamiltonian systems.  For the later we describe the integrals for  the Henon-Heiles Hamiltonian which arises in celestial mechanics.  The unifying feature is that these problems seek elements in the kernel of a linear operator.  The spectral sequence approach emphasizes the obstructions constructed from  cokernel  of the operator   to  finding elements in the kernel.  
	
\end{abstract}

\maketitle

\section{Introduction}\label{sec:section1}

In \cite{bec} a spectral sequence was introduced to study the problem of computing  normal forms for the elements of a graded Lie algebra.  The approach offered a fairly efficient method for specifying  these \color{black} forms in a number of classical 
 contexts, \cite{bec1}, \cite{becc}, each of which could be formulated in terms of the kernel of an adjoint representation. The spectral sequence provided a systematic investigation of this kernel through a chain of approximations.

 In this paper we indicate how those methods can be extended beyond normal forms.  Specifically, we are concerned with uncovering elements of the kernel of a linear mapping $\,f:M \rightarrow N\,$ between modules over a commutative  ring.  The key conceptual results are elementary: the work comes in applications when attempting to verify the hypotheses.The calculations involved in the examples can be extensive, and 
 one can view the role of the spectral sequence as a means for organizing the work.  In particular,
 we will see that certain modules constructed from the cokernel of the relevant linear mapping  can be viewed as obstructions to continuing the calculations. 
  The spectral sequence approach can then be applied to such problems as constructing commuting polynomial vector fields, and determining independent integrals of Hamiltonian vector fields.

The method involves associating one of two spectral sequences with $\,f$,\, and these sequences will come as no surprise to those familiar with such entities.  However, to make the results accessible to a wider audience we do not assume familiarity with such sequences.   In particular we point out that many of  the  notions associated  to a spectral sequences appear with different names in other context.  We illustrate this point  with examples in appendix \eqref{app:ss}  of  problems such  as characterizing  a  Gr\"{o}bner basis for a polynomial ideal and solving polynomial differential equations.

Throughout the paper rings are assumed both commutative and unitary, elements of rings are termed ``scalars."   ``Module'' always means ``unitary left module," and explicit reference to the underlying ring is sometimes omitted.

Our notation and terminology for spectral sequences is essentially that found in \cite[Chapter XI, \S 3 and \S 8]{keymac}.

\section*{Acknowledgment} 

A portion of the work of the third author was done while visiting the University of Calgary.

\section{Vector Fields}\label{section:vectorfields}

\

An application of spectral sequences  is given by 
the graded Lie algebra $\,\cl\,$  of  vector fields over $\mathbf{k}  =\,\nr^n\,$  or $\,\nc^n\,$ having formal power series as component functions and admitting an equilibrium point at the origin.
In dealing with this example we identify a vector field $\,F\,$ with the formally first-order differential equation 
\begin{equation} \label{theeq}
x\pr = F(x).
\end{equation}

We express such an equation 
in terms of column vectors.  That is, if $\,F(x) = (f_1,(x),f_2(x),\dots,f_n(x))\,$ we  write
\begin{equation} \label{colvect}
F = \left[ \begin{array}{c} f_1(x) \\ f _2(x) \\ \vdots \\ f_n(x) \end{array} \right].
\end{equation}

The grading  on $\,\cl\,$ arises by letting   let $\,\cl_j\,$ consist of column vectors of {\it homogeneous} polynomials of degree $\,j+1$.  (The $\,1\,$ within $\,j + 1\,$ is needed to ensure the grading condition in\, (\ref{gradcond}).)

The bracket is the usual Lie bracket.
If a second  vector field $\,G\,$ is identified with $\,x\pr = G(x)$,\, the Lie bracket $\,[F,G]\,$ is identified  with the differential equation
$\,x\pr = E(x)$,\, where
\begin{equation} \label{lbconv}
E(x) = \frac{dG}{dx}(x)F(x) - \frac{dF}{dx}(x)G(x),
\end{equation}
wherein $\,\frac{dG}{dx}(x)\,$ and $\,\frac{dF}{dx}(x)\,$ are the associated Jacobian matrices.  In\, (\ref{theeq})\, and\, (\ref{lbconv})\,   $\,x\,$ denotes the column $[x_1,\dots,x_n]^{\tau}$;\,  when $\,n = 2\,$ we will use $\,(x,y)\,$ in place of $\,(x_1,x_2$)

A classical problem is to characterizing polynomial vector fields that commute with a given polynomial vector field.  In other words for a given vector field $X$  we wish to find vector fields $Y$ such that $[X,Y]=0$.

We apply the methods of Appendix \ref{app:ss} to two examples of vector fields
$F_0 + \cdots + F_p \in \cl $.  We first consider the spectral sequence associated to a decreasing filtration to study vector fields with 

\begin{equation*}
\tag{i} F_0(x) = \left[ \begin{array}{c} x_1 \\ x_2 \\ \vdots \\ x_n \end{array} \right]
\end{equation*}

\nd we then consider the spectral sequence associated to an increasing filtration to study vector fields  with

\begin{equation*}
\tag{ii} F_p(x) = \dfrac{1}{p+1}\left[ \begin{array}{c} x_1^{p+1} \\ x_2^{p+1} \\ \vdots \\ x_n^{p+1} \end{array} \right].
\end{equation*}

\nd The calculations are based on:

\begin{lem}[Euler's Formula]\label{euler} {\rm

	Suppose $F:\nr^n \to \nr^m$ is a differentiable homogeneous function of degree $r$ one has

	\begin{equation*} \tag{i}
		\dfrac{dF}{dx} \cdot x = r F(x) \mbox{ for all } x \in \nr^n.
	\end{equation*}

	where 	
	
	\begin{equation*}\tag{ii}
		\dfrac{dF}{dx}(x) \cdot x \eqqcolon \Big[\dfrac{df_i}{dx_j}\Big]\left[ \begin{array}{c} x_1 \\ x_2 \\ \vdots \\ x_n \end{array} \right].
	\end{equation*}
	
}
	
\end{lem}

\

\pf By assumption one has 

\begin{equation*}
	\tag{iii} F(\lambda x)= \lambda^r F(x).\end{equation*}

By differentiating with respect to $\lambda$ we have:

$$\dfrac{d F(\lambda x)}{d\lambda}  =  \dfrac{dF}{dx} (\lambda x) \frac{d}{d\lambda}\lambda x   = \dfrac{dF}{dx}(\lambda x) \cdot x =   r \lambda ^{r-1} F(x),$$

\nd whereupon (i) is obtained by setting $\lambda =1$. \qed

\begin{cor}\label{prop:ad}{\rm
	Suppose $t,s \in \nn$  and $F \in \cl_t, G \in \cl_s$.  Then \medskip

$$	{\rm ad}(F)(G) =\displaystyle{}\frac{1}{1+t}\,\Big( \dfrac{dG}{dx}\, \dfrac{dF}{dG} - \displaystyle{}\ell\, \dfrac{dF}{dX}\, \dfrac{dG}{dx}\Big) \cdot x,  \quad \text{where} \quad \ell := \frac{1+t}{1+s}.$$}
\end{cor}
\proof

 From Lemma \ref{euler} one has\medskip
\begin{equation*}\tag{i}F= \dfrac{1}{t+1} \cdot \dfrac{dF}{dx} \cdot x \mbox{ and }
G= \dfrac{1}{s+1} \cdot \dfrac{dG}{dx} \cdot x. \medskip 
\end{equation*}
The corollary now follow immediately from (i) and the equalities below:
	\begin{equation*}
		\begin{array}{lllllll}
			{\rm ad}(F)(G) & = & \dfrac{dG}{dx} F - \dfrac{dF}{dx} G \medskip \\
			& = & \displaystyle{}\frac{1}{1+t}\, \dfrac{dG}{dx}\, \dfrac{dF}{dG} \cdot x - \displaystyle{}\frac{1}{1+s}\, \dfrac{dF}{dX}\, \dfrac{dG}{dx} \cdot x\medskip \\
			& = & \displaystyle{}\frac{1}{1+t}\,\Big( \dfrac{dG}{dx}\, \dfrac{dF}{dG} \cdot x - \displaystyle{}\ell\, \dfrac{dF}{dX}\, \dfrac{dG}{dx} \cdot x\Big), \quad \text{where} \quad \ell := \frac{1+t}{1+s}, \medskip \\
			& = & \displaystyle{}\frac{1}{1+t}\,\Big( \dfrac{dG}{dx}\, \dfrac{dF}{dG} - \displaystyle{}\ell\, \dfrac{dF}{dX}\, \dfrac{dG}{dx}\Big) \cdot x,  \quad \text{where} \quad \ell := \frac{1+t}{1+s}.
		\end{array}
	\end{equation*} \smallskip
\qed

\begin{cor}\label{cor46}{\rm
The following statements are equivalent.

\begin{itemize}
	\item[(a)] $({\rm ad}\, F)(G) = 0$, i.e. $\,G \in \ker ({\rm ad}\,F)$.\medskip
	\item[(b)] $\dfrac{dG}{dx} F = \dfrac{dF}{dx} G.$ \medskip
	\item[(c)] $\big(\dfrac{dG}{dx} \dfrac{dF}{dx} - \ell \dfrac{dF}{dx}\dfrac{dG}{dx}\big) \cdot x = 0$, where
	$\ell := \displaystyle{} \frac{t+1}{s+1}$. \medskip
	\item[(d)] $M(x) \cdot x = 0$,\, where $\,M = M(x) = \dfrac{dG}{dx}(x) \dfrac{dF}{dx}(x) - \ell \dfrac{dF}{dx}(x)\dfrac{dG}{dx}(x)$. \medskip
	\item[(e)] 	If $\ \dfrac{dF}{dx}\,$ is invertible, an additional equivalent condition is \medskip
	
		  $ \left( \left(\dfrac{dF}{dx}\right)^{-1}\,\dfrac{dG}{dx}\,\dfrac{dF}{dx}\right) \cdot x = \ell \, \dfrac{dG}{dx} \cdot x. $
\end{itemize}
}
\end{cor}

\qed
\medskip

\nd Note that the matrix $\,M\,$ appearing in\, (d)\, is either $0\,$ or homogeneous of degree $\,s+t.$

\begin{prop}\label{prop:adiso}{\rm
	Suppose
	 \begin{equation*}
	 F_0(x) = \left[x_1,x_2,\cdots,x_n\right]^{\tau}.
	\end{equation*}Then

	\label{prop:inj}
	$$ad(F_0)|_{\cl_j}: \cl_j \to \cl_j$$

	\nd is injective for $j \geq 1$ and $ad(F_0)(\cl_0) = 0$.
}
\end{prop}

\pf We have $\dfrac{dF_0}{dx} = I$.  For any $0 \neq G\in \cl_j$ we have by Corollary \ref{cor46}

\begin{equation}\label{adjointmap}ad(F_0)(G) = \dfrac{dG}{dx}\cdot x -\dfrac{1}{j+1}\dfrac{dG}{dx}\cdot x=(\frac{j}{j+1} )\dfrac{dG}{dx} \cdot x = \frac{j^2}{j+1} G.\end{equation}

\qed

The graded  Lie algebra $\cl = \underset{j\geq 0}{\prod} \cl_j$ induces   a spectral sequence,  $E_r^{p,s}$ (See \eqref{def:decss}). Feeding Proposition \ref{prop:inj} into the spectral sequence enables us to characterize all centralizers of a polynomial vector field with lowest degree term $[x_1,x_2, \cdots, x_n]^{\tau}$.

\begin{cor}{\rm
	\label{cor:formalintegral} Suppose $F=F_0+F_1 + \cdots $ where
	 \begin{equation*}
	 F_0(x) = \left[x_1,x_2,\cdots,x_n\right]^{\tau},
	\end{equation*}

	 \nd then	every $G_0 \in \cl_0$ extends to a formal centralizer of $F$:

	$$G = G_0 + G_1 + \cdots, \quad ad(F)(G)=0, $$ 
	
	\nd where $G_i$ is inductively given by:
	
	\begin{equation}\label{central}G_i = -\dfrac{i+1}{i^2}\Big( [ F_1, G_{i-1}] + [ F_2, G_{i-2}] + \cdots + [ F_i, G_{0}] \Big).\end{equation}
}

\end{cor}

\pf  Since $\cl_j$ is a finite vector space for all $j$ Proposition \ref{prop:inj} implies  $d_0:E_0^{p,0} \to E_0^{p,1} $ is an isomorphism for all $p \geq 1$.   Proposition \ref{viewtwo}  implies $G$ exist.

To prove \eqref{central} We inductively assume that $[G_0] \in  E_0^{0,0}$ survives to $E_{r-1}^{0,0}$.  This means that $G_j \in \mathcal{L}_j$ have been defined for $0\leq j <k$ so that the filtration of $ad(F)(\underset{i=0}{\overset{k-1}{\Sigma}} G_i)$ is greater than $k-1$.

In other words
$$d_{r-1}([G_0]) = ad(F)(\sum_{i=0}^{k-1} G_i) = \sum_{i=0}^{k-1} [F_{k-i},G_i] +\mbox{ plus terms of higher degree. } $$

By Proposition \ref{prop:adiso} \hspace{.1 in} $d_{r-1}$ takes values in the $0$ group. Specifically from \eqref{adjointmap}  \hspace{.1 in} $\sum_{i=0}^{k-1} [F_{k-i},G_i] \in \mathcal{L}_k$ is   $ad(F_0)$ of the element

$$  \dfrac{k+1}{k^2} \sum_{i=0}^{k-1} [F_{k-i},G_i]. $$

Defining $G_k$ by \eqref{central} completes the induction.

\qed

\begin{exm}{\rm 
	
	Clearly there are finite centralizers of $F$, namely $F$ itself.  In this example we show that there are formal infinite centralizers.  Let $F$ be the vector field
	
	$$ F_0+F_1=
	 \left[x_1,x_2,\cdots,x_n\right]^{\tau} + \dfrac{1}{2}\left[x_1^2,x_2^2,\cdots,x_n^2\right]^{\tau}.
	$$

By Corollary \ref{cor:formalintegral} each $G_0 \in \mathcal{L}_0$ extends to a centralizer with each $G_i$ an iterated bracket
$$[F_1,[F_1 \cdots [F_1,G_0] \cdots ].$$ up to a no-zero factor.

We will prove in Proposition \ref{prop:injgr} that $ad(F_1):\mathcal{L}_t \to \mathcal{L}_{t+1}$ is injective for $t\neq 1$.  As a consequence any 
centralizer with $G_2 \neq 0$ must be an infinite sum.  For example for $n=2$ and $G_0=[y, x]^{\tau}$ one has $G_2 \neq 0$.

}
\end{exm}

For our  next examples we study formal vector fields $F=F_0+ \cdots F_p$ where 
$F_p= \dfrac{1}{p+1} [x_1^{p+1},\cdots, x_n^{p+1}]^{\tau}$.  In the sequel $F$ is such a formal vector field.

Our main result is that any finite centralizer of such a vector field must be a vector field of polynomials of degree $p$.

\begin{thm}\label{thm:adFr}{\rm  
	
	There does not exist a centralizer, $G$ of $ad(F)$ of the form
	$$ G = G_0 + \cdots G_j, \quad j\neq p, \quad G_j \neq 0$$}
\end{thm}

Theorem \ref{thm:adFr} is a consequence of the following proposition.

\begin{prop}\label{prop:injgr}{\rm   Suppose $F=F_0 + \cdots + F_p \in \cl $ with $F_p =  \dfrac{1}{p+1} [x_1^{p+1},\cdots, x_n^{p+1}]^{\tau}$.  Then

			$$ad(F_p)|_{\cl_j}: \cl_{j} \to \cl_{p+j}$$ is injective if $j \neq p$.

}
	
\end{prop}

\nd \textbf{Proof of Theorem \ref{thm:adFr} (assuming Proposition \ref{prop:injgr}): }  Let $\widehat{\cl}$ be the graded $\mathbf{k}$ module with 

$$\widehat{\cl}_j = \cl_{j-p}.$$  In particular $\widehat{\cl}_j = 0$ if $j<p$ and  $\widehat{\cl}_p = \cl_0$.  $\widehat{\cl} = \underset{j \geq 0}{\prod} \widehat{\cl}_j$ is a graded $\mathbf{k}$, module  and 
$$f=ad(F_p): \widehat{\cl} \to \cl$$ is a map of graded modules.

We consider the spectral sequence $E^r_{p,*}$ associated to an {\it increasing} filtration of $f$ mentioned in   \eqref{def:ssinc} of Appendix \ref{app:ss}. Proposition \ref{prop:injgr} implies $d_0:\widehat{\cl}_j \to \cl_j$ is injective for $j \neq 2p$.  The Theorem now follows from Corollary \ref{keycor}.

\qed

\nd The following lemma is used to prove 
Proposition \ref{prop:injgr}.

\begin{lem}\label{lem:312}{\rm 
	
	Suppose $g(x)$ is a homogeneous polynomial of degree $>0$ in $n-$variables which satisfies

	$$ \sum_{j=1}^n \dfrac{\partial g}{\partial x_i}x_i^{p+1} = 0.$$

\nd	Then $g(x) =0$.
}
	
\end{lem}

\pf  If $p=0$ then $\sum_{j=1}^n \dfrac{\partial g}{\partial x_i}x_i = s g$ where $s$ is the degree of homogeneity. So we may assume $p>0$.

Assume to the contrary that $g(x) \neq0$.  Then there is a variable $x_i$ that appears in $g$.  Without loss of generality we may assume $i=1$.  We can write

$$g(x) = x_1^r f_r(x_2, \cdots, x_n) + x_1^{r-1} f_{r-1}(x_2, \cdots, x_n) + \cdots + f_0(x_2, \cdots, x_n)$$\\

with $r>0, \quad f_r \neq 0$.  Then

$$ r x_1^{r+p} f_r + (r-1) x_1^{r+p-1} f_{r-1} + \cdots x_1^{p+1}f_1 
+ \sum_{i=2}^n \dfrac{\partial g}{dx_i} x_i ^{p+1}=0.$$

Since $p>0$, the largest power of $x_1$ appears only in the first term, hence cannot be canceled, thereby implying the contradiction $f_r = 0$.

\qed

\begin{lem} \label{matone}{ \rm Suppose $\,\dfrac{dF}{dx}\,$ is diagonal, say
	\begin{equation*} \label{i}
	\dfrac{dF}{dx} = {\rm diag}(\lambda_1,\dots,\lambda_n) =: \lambda. 
	\end{equation*}

\nd Write $G$ as $[g_1,\cdots, g_n]^{\tau}$. Then the column vector

	$$\Big[\left( \dfrac{dG}{dx} \dfrac{dF}{dx} - \ell\,\dfrac{dF}{dx}\, \dfrac{dG}{dx}\right) \cdot x\Big]$$ has i-th entry \medskip
	\begin{equation*} 
	 \begin{array}{lllll}
	\displaystyle{\sum_j(\lambda_j-\ell \lambda_i) \frac{\partial g_i}{\partial x_j} x_j }.
	\end{array}
	\end{equation*}
\nd where $\ell$ is defined in Corollary \ref{prop:ad}.
}
\end{lem}

\pf Write $\,F\,$ as $\,[f_1,f_2,\dots,f_n]^{\tau}$.   Since $\,\dfrac{dF}{dx}\,$ is diagonal each function, $f_i(x)$ is a function of $x_i$.
 To ease notation we write $\,\dfrac{dG}{dx}\,$ as $\,[g_{ij}]\,$ and 
$\dfrac{dF}{dx} = [\delta_{ij}\,\lambda_j]$.

\nd We have 
\begin{equation*}
\big(\dfrac{dG}{dx} \dfrac{dF}{dx}\big)_{ij} = \sum_k g_{ik}\delta_{kj}\lambda_j = \lambda_jg_{ij},
\end{equation*}

\begin{equation*}
\big(\dfrac{dF}{dx} \frac{dG}{dx}\big)_{ij} = \sum_k\delta_{ik}\lambda_i g_{kj} = \lambda_ig_{ij}.
\end{equation*}

\nd and
\begin{equation*} \tag{i}\label{tagi}
\dfrac{dG}{dx} \dfrac{dF}{dx} - \ell\,\dfrac{dF}{dx} \dfrac{dG}{dx} = [(\lambda_j-\ell\lambda_i)\frac{\partial g_i}{\partial x_j}].
\end{equation*}

\nd Hence the $\,i^{\rm th}$-entry of the column vector $\, \left( \dfrac{dG}{dx} \dfrac{dF}{dx} - \ell\,\dfrac{dF}{dx}\, \dfrac{dG}{dx}\right) \cdot x\,$ is
\begin{equation*}
\displaystyle{\sum_j(\lambda_j-\ell \lambda_i) \frac{\partial g_i}{\partial x_j} x_j },
\end{equation*}
which establishes\, the lemma.
\qed

\begin{thm} \label{thm:416}{\rm  Assume $F_p=  \dfrac{1}{p+1} [x_1^{p+1},\cdots, x_n^{p+1}]^{\tau}
	\in \mathcal{L}_p, \,\newline p>0$.

	\begin{enumerate}
	
\item If $\ell \neq 1$	
	then the mapping $\,{\rm ad}(F):L_t \rightarrow L_{p+t}\,$ is injective.

\item $\text{ker}(F_p):\mathcal{L}_p \to \mathcal{L}_{2p}$ is generated by
$$g_i = [0,\cdots, x^{p+1}, \cdots,0]^{\tau}, \mbox{  where } x^{p+1} \mbox{   is in the i-th position}$$	
	\end{enumerate}
}
\end{thm}

  \bigskip

\pf We write $\,
\dfrac{dF}{dx} = {\rm diag}(x^p_1,\dots,x^p_n)\,$.

Assume  $G \in \mathcal{L}_t, t \neq p$ is in the kernel of $\,{\rm ad}\,(F)$. Then by\, Corollary \ref{cor46} and  Proposition \ref{matone}\, we have

\begin{equation*}
\displaystyle{\sum_j(x^p_j-\ell x^p_i) \frac{\partial g_i}{\partial x_j} x_j } =0, \, i=1, \cdots, n.
\end{equation*}
Since each $\,g_i\,$ is homogeneous of degree $\,t+1\,$ we can appeal to Lemma \ref{euler} to express this as
\begin{equation*} \tag{ii}
(p+1)\,x_i^p\, g_i =  \sum_j x_j^{p+1}\frac{\partial g_i}{\partial  x_j}  \qquad \text{for all indices $\,i$.}
\end{equation*}

If $\,G \neq 0\,$ at least one $\,g_i\,$ must be non-zero, and some $\,x_j\,$ must appear in $\,g_i$, w.l.o.g. $\,j = 1$.  Express $\,g_i\,$ as a polynomial in $\,x_1\,$ with coefficients in $\,A[x_2,\dots,x_n]$, say
\begin{equation*}
\hspace{-1.3 in} g_i = x^r_1h_r(x_2,\dots,x_n) + x^{r-1}_1h_{r-1}(x_2,\dots,x_n) + $$$$\cdots + x_1h_1(x_2,\dots,x_n) + h_0(x_2,\dots,x_n),
\end{equation*}
where $\,h_r \neq 0.$

 We can then write (ii) for this particular index $i$, as

\begin{equation*} \tag{iii}\label{tagiiiprop}
(p+1)   x_i^p [ x^r_1 h_r + \cdots ]= x_1^{p+1} [rx_1^{r-1}h_r + \cdots]+x_2^{p+1} [x_1^r \frac{\partial h_{r-1}}{\partial x_2} + \cdots]+\cdots
\end{equation*}

\nd Suppose $i \neq 1$.  Since $\,h_r \neq 0\,$ we see  that the highest power of $\,x_1\,$ on the left in\, \eqref{tagiiiprop}\, is $\,r$,\,  whereas that on the right  is $\,r+p$.  Since $p>0$ this forces the contradiction $\,h_r = 0$.

\nd We now assume $i=1$ then \eqref{tagiiiprop} becomes

\begin{equation*} \tag{iv}\label{tagivprop1}
(p+1)   x_1^{p} [ x^r_1 h_r + \cdots ]=   x_1^{p+1}[rx_1^{r-1}h_r + \cdots]+x_2^{p+1}  [x_1^r \frac{\partial h_{r}}{\partial x_2} + \cdots]+\cdots
\end{equation*}

The highest power of $x_1$ on both sides is $p+r$, with coefficient $(p+1) h_r $ on the left and coefficient $ r h_r$ on the right.  So if  $r \neq p+1$ we again have the contradiction $h_r =0$.

Now assume $r=p+1$ and, after  rearranging the terms in \eqref{tagivprop1} we have

\begin{multline*} \tag{iv}\label{tagivprop2}
(p+1)   x_1^{p} [ x^{p+1}_1 h_{p+1} + \cdots ]-  x_1^{p+1}[(p+1)x_1^{p}h_{p+1} + \cdots]= \\+ x_2^{p+1} [x_1^{p+1} \frac{\partial h_{p+1}}{\partial x_2} + \cdots]+x_3^{p+1}[x_1^{p+1} \frac{\partial h_{p+1}}{\partial x_3 \cdots} + \cdots]\cdots .\\
\end{multline*}

The highest power of $x_1$ on the right is $p+1$  so all the terms with higher powers of $x_1$ on the left must be zero.  These terms are:

$$\sum_{k=2}^{p} (p+1) x_1^{p+k}h_k- x_1^{p+1} (k x_1^{k-1}h_k)=\sum_{k=2}^p (p-k+1)h_k x_1^{p+k},   $$

\nd which implies $h_k=0, \quad 2 \leq k \leq p$. So $g_1$ has the form

\begin{equation*}\tag{v}\label{prooftagv} g_1= x_1^{p+1} h_{p+1} + x_1 h_1 + h_0.\end{equation*}

\nd Substituting \eqref{prooftagv} into the left side of \eqref{tagivprop2} gives:

$$(p+1)\Big[x_1^{2p+1}h_{p+1}+x_1^{p+1}h_1+x_1^p h_0 \Big] - \Big[ (p+1) x_1^{2p+1}h_{p+1}+x_1^{p+1}h_1\Big]$$ $$
=(p+1)x_1^{p+1}h_1+(p+1)x_1^ph_0-x_1^{p+1}h_1
$$$$
=px_1^{p+1}h_1+(p+1)x_1^ph_0.$$

\nd Substituting \eqref{prooftagv} into the right side of \eqref{tagivprop2} gives

$$x_1^{p+1}(\sum_{i>1} \dfrac{\partial h_{p+1}}{\partial x_i} x_i^{p+1})+x_1(\sum_{i>1} \dfrac{\partial h_1}{\partial x_i} x_i^{p+1})+  (\sum_{i>1} \dfrac{\partial h_0}{\partial x_i} x_i^{p+1}). $$

Since $p$ is the lowest power of $x_1$ on the left side we must have 

$$ \sum_{i>1} \dfrac{\partial h_0}{\partial x_i} x_i^{p+1}= \sum_{i>1} \dfrac{\partial h_1}{\partial x_i}x_i^{p+1} =0.$$

\nd Lemma \ref{lem:312} implies $h_0=h_1=0$.

\nd The left side  is now zero, which implies

$$\sum_{i>1} \dfrac{\partial h_{p+1}}{\partial x_i} x_i^{p+1}=0$$

\nd Since the homogeneity degree of $h_{p+1} $  is strictly grater than zero, Lemma \ref{lem:312}  implies $h_{p+1}=0$ which is a contradiction \footnote{This is the only place we are using $\ell \neq 1$.}  proving (1).

If $t=p$ then (2) follows since the only element in the kernel is given by $h_{p+1}$ a constant.

\qed

There is the spectral sequence associated to the increasing filtration of $\mathcal{L}$ described in Appendix \ref{app:ss}.  Briefly the elements in $F^j\mathcal{L}$ of filtration $j$ are the classes, $m \in \mathcal{L}$ of the form
$$m=m_0+\cdots +m_j.$$

For $F=F_0+ \cdots F_p \in \mathcal{L}$ \, there is the linear map $ad(F):F^j\mathcal{L} \to F^{p+j}\mathcal{L}$.  The map $ad(F)$ becomes a graded linear map if we re-index $\mathcal{L}$ by 
$\widehat{\mathcal{L}}_j=\mathcal{L}_{j-p}$ \,
as in the proof of Theorem \ref{thm:adFr}.

\nd The  spectral sequence associated to the graded linear map
$$\ad(F):\widehat{\mathcal{L} }\to \mathcal{L}$$ is denoted $E^r_{p,q}$.  For $F=F_0 + \cdots + F_p \in \cl $ with 

\begin{equation*}
 F_p(x) = \dfrac{1}{p+1}\left[  x_1^{p+1}, x_2^{p+1}, \cdots , x_n^{p+1}  \right]^{\tau}.
\end{equation*}

\nd Proposition \ref{prop:injgr} implies 
$E^1_{s,-s} =0 \mbox{ if } s \neq 2p.$

We now specialize to 
\begin{equation}\label{F0Fp} F=F_0+F_p \mbox{ where } F_0=[x_n, x_{n-1}, \cdots,x_1]^{\tau}.\end{equation} 
 The differentials in the spectral sequence for this vector field will be determined and, as a consequence we will have computed the vector space of centralizers of $F$. In the following theorem the $f_i$ denotes the coordinates of  functions $F$, e.g. $F=[f_1,\cdots,f_n]^{\tau}$.

  \begin{thm}\label{thm:basisforker}  With $F$ as in \eqref{F0Fp}, 
 	$\text{ker}(ad(F))$ has as a basis the classes

 	$$[0,\cdots, f_i, 0,\cdots,0, f_{n-i+1}, \cdots, 0]^{\tau}, \quad 1 \leq i \leq \frac{n+1}{2}.$$
 	
 \end{thm}

This is a generalization of the example found in \cite{not} where they study the vector field

$$\left[ \begin{matrix}
y+ x^3 \\
x+ y^3
\end{matrix}\right].$$

To prove Theorem \ref{thm:basisforker} we start with the basis for $E_{2p,0}^1 \subset \widehat{\cl}_{2p} = \cl_p$ given by Theorem \ref{thm:416}.

\begin{equation}\label{def:gi} \{ \alpha_i:=\left[  0 , \cdots, 0 , x_i^{p+1}, 0,\cdots,0  \right]^{\tau}  \quad i=1, \cdots, n\}.\end{equation}

\nd where $ x_i^{p+1}$ is in the $i-$th row.

Suppose $G=[g_1,\cdots, g_n]$ is a finite centralizer of $F$.  An immediate consequence of  \eqref{def:gi}  is that

 $$g_i=\sum_{j=1}^n c_i x_i^{p+1} \mbox{ modulo lower degree}.$$

The {\it lower degree terms } are determined by the first non-trivial differentials of $ \alpha_i $ which is given by

\begin{equation} \label{def:dp} d_p(\alpha_i)\footnote{It is a $d_p$ since it drops filtration by p.}   = ad(F_0)(\alpha_i) + ad(F_p)(\alpha_i)= ad(F_0)(\alpha_i).\end{equation}

  which takes values in $$E^p_{p,1} = \cl_p/\text{im}(d_0) =\cl_p/(\text{im}(ad(F_p)|_{\widehat{\cl}_p}).$$

To compute $\cl_p/(\text{im}(ad(F_p)|_{\widehat{\cl}_p})$ we evaluate 
$ad(F_p)$ on a basis, $h_{i,j}$,  for $\cl_p$, obtaining

$$ h_{i,j} = \left[  0 ,\cdots, 0,  x_i, 0, \cdots, 0  \right]^{\tau}  \quad \mbox{ with } x_i \mbox{ in the } j-th \mbox{ row }.$$

\begin{lem}\label{lem:adhij}
    $$(p+1)ad(F_p)(h_{i,j}) = x_i^ph_{i,j} - (p+1) x_j^p h_{i,j}.$$
\end{lem}

\qed

\begin{cor}\label{cor:rel}

$$E^p_{p,1} =\cl_p/I,$$

\nd where $I$ is the ideal generated by the following relations in the $j-$th row: 

$$ \left[ \begin{array}{c} 0 \\ \vdots\\ 0 \\ x_i^{p+1} \\0\\\vdots\\0 \end{array} \right] =  \left[ \begin{array}{c} 0 \\ \vdots\\ 0 \\ (p+1)x_j^px_i \\0\\\vdots\\0 \end{array} \right]  \quad i =1, \cdots, n.$$
\end{cor}

\qed

\nd As an example, suppose  $i=j$ and  $x_i^{p+1}=0$, in the $i-$th row.  Corollary \ref{cor:rel} implies that  $E^p_{p,1}$ is generated in each row by 

$$\{ x_1^{\epsilon_1}\cdots x_n^{\epsilon_n} \Big| 0 \leq \epsilon_i \leq p\}.$$

We now determine the kernel of  $d_p $.  To this end we have the following lemma:

\begin{lem}\label{lem:dpgi}{\rm 
	
	$$(p+1)d_p(\alpha_i) =  \left[ \begin{array}{c} 0 \\ \vdots\\ 0 \\ (p+1)x_i^px_{n-i+1} \\0\\\vdots\\-x_i^{p+1}\\\vdots\\0 \end{array} \right]= \left[ \begin{array}{c} 0 \\ \vdots\\ 0 \\ (p+1)x_i^px_{n-i+1} \\0\\\vdots\\-(p+1)x_{n-i+1}^p x_i\\\vdots\\0 \end{array} \right], $$
	
	\nd where $x_i^px_{n-i+1}$ is in the $i-$th row and $x_i^{p+1}$
	is in the $(n-i+1)-$st row.
}
\end{lem}

\pf The calculation of the differential is immediate from the formula for $ad(F_0)(\alpha_i)$.  The second equality follows from the relations in Corollary \ref{cor:rel}. \qed

\begin{cor}\label{cor:keradF}{\rm 
	$E^{p+1}_{2p,0}$ has a basis the classes  $$\{ \alpha_i+ \alpha_{n-i+1}, \quad 1 \leq i \leq \frac{n+1}{2} \}.$$
}
\end{cor}

\pf By Lemma  \ref{lem:dpgi} the classes $\alpha_i+ \alpha_{n-i+1} $ generate the kernel of $d_p$. 

\qed

\nd  We next outline the computation of the next (and final) differential.  We need representatives for $\alpha_i+\alpha_{n-i+1}$ in $E^{p+1}_{2p,0}$.   We have

 \begin{multline*}\hspace{-.2 in}  ad(F_0)(\alpha_i + \alpha_{n-i+1})=\\
 [0,\cdots,(p+1)x_i^px_{n-i+1}-x^{p+1}_{n-i+1},0 \cdots,-x_i^{p+1} +(p+1)x^p_{n-i+1}x_i,0,\cdots,0]^{\tau}.  \\ \end{multline*}

\nd This class is zero in $$E_{p,1}^p =\cl_p/(\text{im}(ad(F_p)|_{\widehat{\cl}_p}.$$
 
\nd It follows that  there is a class $\beta \in \widehat{\cl}_p$ with

 $$\hspace{-5.2 in} ad(F_p)(\beta)=$$$$[0,\cdots,(p+1)x_i^px_{n-i+1}-x^{p+1}_{n-i+1},0 \cdots,-x_i^{p+1} +(p+1)x^p_{n-i+1}x_i,0,\cdots,0]^{\tau}.$$

\nd The class $\beta$ is indeed the class that generated the relation in Lemma \ref{lem:dpgi}:
 
 $$-ad(F_p)(h_{n-i+1,i}+ h_{i, n-i+1}) = d_p(\alpha_i+\alpha_{n-i+1}).  $$

 In other words, the class

 $$\alpha_i+\alpha_{n-i+1} + h_{n-i+1,i}+ h_{i, n-i+1} \in \cl$$ satisfies
 \begin{itemize}
 	\item $[\alpha_i+\alpha_{n-i+1}] = [\alpha_i+\alpha_{n-i+1} + h_{n-i+1,i}+ h_{i, n-i+1}] \in E^{p+1}_{2p,-2p}$

 	\nd because $ h_{n-i+1,i}+ h_{i, n-i+1}$ has lower degree than $\alpha_i+\alpha_{n-i+1}$.

 	\item The class has been chosen so that   $ ad(F)(\alpha_i+\alpha_{n-i+1} + h_{n-i+1,i}+ h_{i, n-i+1})$ has filtration $p-1$.  In fact $ ad(F)(\alpha_i+\alpha_{n-i+1} + h_{n-i+1,i}+ h_{i, n-i+1})=0$.

 	\item $$\hspace{-1.8 in}\alpha_i+\alpha_{n-i+1} + h_{n-i+1,i}+ h_{i, n-i+1} = $$$$
 	[0,\cdots, f_i, 0,\cdots,0, f_{n-i+1}, \cdots, 0]^{\tau}$$
 
 	\end{itemize}

 This proves Theorem \ref{thm:basisforker}.
 
 \qed
 
 \section{Application to Hamiltonian Vector Fields}\label{sec:appHam}

 See appendix \ref{sec:hamiltonian} for the notation and a  brief introduction of Hamiltonian vector fields. In particular $\ch$ denotes the Lie algebra of formal power series in the variables $$x=(x_1, \cdots,x_{2n})=(q_1,q_2,\cdots.q_n,p_1,p_2,\cdots,p_n)$$ with Lie bracket given by the Poisson bracket (defined  in \eqref{def:poisson}).   There is a grading on $\ch$ with $\ch_i$ equal to the homogeneous monomials of degree $i+2$.

 We first observe that Hamiltonian vector fields are not an example of the vector fields studied in Section \ref{section:vectorfields}.

 \begin{prop} \label{wontwork}{\rm The hypotheses of\, {\rm Theorem \ref{matone} } fail when $\,t \geq 2\,$ and $\ F:R^{2n} \rightarrow R^{2n}\,$ is a Hamiltonian system.  
 	}
 \end{prop}
 \proof  If $\,\dfrac{dF}{dx} = \lambda\,$ then $\lambda_i=a_i x_i^p$ for $a_i \in A$.  Since $\,F\,$ is a Hamiltonian system with Hamiltonian $\,H\,$ we have $\,F = JH_x$ and
 \begin{equation*}
 \lambda = \dfrac{dF}{dx} = JH_{xx}.
 \end{equation*}
 Multiplying by $\,-J = J^{-1}\,$  then gives
 \begin{equation*}
 -J\lambda = H_{xx}.
 \end{equation*}
 But $\,H_{xx}\,$ is symmetric, whereas $\,-J\lambda\,$ is symmetric if and only if
 \begin{equation*} \tag{i}
 a_jx_j^p = -\lambda_{n+j} \qquad \text{for} \qquad j = 1,\dots,n.
 \end{equation*}
 However, for $\,t \geq 2\,$ the only variable appearing in $\,x_j^p\,$ is $\,q_j$,\, and the only variable appearing in $\,\lambda_{n+j}\,$ is $\,p_j$. Equality\, (i)\, is therefore impossible. \qed

We   study centralizers of the Henon-Heiles Hamiltonian which have non-zero, quadratic lowest terms. 

 \begin{defin}\label{def:henon-heiles}{\rm 
	The {\it Henon-Heiles} Hamiltonian is the vector field on $\nr^4$ defined by

	$$H=H(q_1,q_2,p_1,p_2) = \dfrac{1}{2}(p_1^2+p_2^2) + \dfrac{1}{2}(A q_1^2 + B q_2^2) + \dfrac{1}{3} q_1^3 + L q_1 q_2^2.$$

}
\end{defin}

In term of the filtration on $\ch$,

$H= H_0 + H_1 \mbox{ where }$$$H_0 = \dfrac{1}{2}|p|^2+ \dfrac{1}{2}(A q_1^2 + B q_2^2) \mbox { and } H_1 = \dfrac{1}{3} q_1^3 + L q_1 q_2^2.$$

 The spectral sequence approach to this problem breaks the computation into a number of steps.

\begin{itemize}
	
	\item View $ad(H):\ch \to \ch$ as a graded linear map with the decreasing filtration defined in \eqref{asgrdcr}.\\

	\item  Compute $E^{0,0}_1 =\text{ker}(ad(H_0)\Big|_{\ch_0}$.  These are the potential quadratic lowest terms of centralizers of $H$.\\

	\item Compute the higher differentials on  $E_r^{0,0}$ for $r\leq 3$. This determines all centralizers of $H$ of degree $\leq 4$ with quadratic lowest term. 
	\end{itemize}  

\nd It is convenient to introduce complex variables which are chosen so that the quadratic part of the Henon-Heiles equation has a nice form: $$H_0= \frac{1}{2} ( z_1 \zbar_1 + z_2 \zbar_2).$$. 

\nd where 
 
\begin{equation}\label{cplxHamgenerators}
\begin{array} {cccc}
z_1 &= & \sqrt{A} q_1 + i p_1  \\
\vspace{.2 in}
\zbar_1 &=& \sqrt{A} q_1 - i p_1\\

z_2 &= & \sqrt{B} q_2 + i p_2  \\
\zbar_2 &=& \sqrt{B} q_2 - i p_2 .

\end{array}\end{equation}

\nd This is a basis for $\nc^4$ if $AB\neq0$, which we now assume.

  We next  describe the homomorphisms
  
  $$ad(H_0):\ch_i \to \ch_i.$$

\nd The kernel of $ad(H_0)$ is given by Corollary \ref{cor:5.10}  if $i=0$.  We need the cokernel of $ad(H_0)$ for $i>0$ in order to compute $E_1^{p,1}$ which is the input for computing the higher differentials in the spectral sequence.

 There is no extra work to compute $ad(H_0)$ for the natural generalization of $H_0$ to   $\nr^{2n}=(q_1,\cdots, q_n,p_1,\cdots,p_n)$:

    $$H_0=\dfrac{1}{2}\sum_{j=1}^n p_j^2 + \dfrac{1}{2} \sum_{j=1}^n A_j q_j^2.$$

  \nd The generalization of  the complex generators  \eqref{cplxHamgenerators}
  is $z_i = \sqrt{A_i} q_i+i p_i, \quad i=1, \cdots n$.

\nd If $n=2$, $A_1=A$ and  $ A_2=B$.

   We need a definition.

 \begin{defin}\label{def:[m]}{\rm 
For a homogeneous monomial 
$m=\prod\limits_{j=1}^n z_j^{I_j}\zbar_j^{J_i} \in \ch_{(\sum I_j +\sum J_j)-2}$  define $[m]\in \nc$ by

$$[m] = \sum_{j=1}^n k_j \sqrt{A_j} ,$$

\nd where $k_j= J_j- I_j.$
}

 \end{defin}

 \begin{thm}\label{thm:conj} {\rm

 	\begin{enumerate}[label=(\roman*)]
 		
 		\item $$ad(H_0)(m) = (i[m]) m.$$

 		\item  The centralizers of $H_0$ are freely generated by monomials $m$ with $[m]=0$.

 		\item The  cokernel of  $ad(H_0)$ is freely generated by monomials, $m$ \newline with $[m]=0$.

 	\item The centralizers of $H_0$ are closed under multiplication.

 	\item	
 	Suppose $A_j$  satisfy the following condition:
 	
 		\begin{equation}\label{hyp:conj}
 \sum_{j=1}^n k_j \sqrt{A_j}  \neq 0 \mbox{ for all  }k_j\in \nz, \mbox { such that  } \sum k_j = 1 \mod 2\end{equation}
 
 Then all centralizers are of  even degree.

 	\item Suppose $n=2$ and $A=B$.  Then $P $ is a polynomial algebra on the generators of $\text{ker}(ad(H_0)\Big|_{\ch_0}$.  In this case $E_1^{2p+1,1}=0$.

 	\end{enumerate}
}
 	
 	\end{thm}

 \begin{rmk}{\rm 
 
 \
 
 \begin{enumerate}
\item  Condition \eqref{hyp:conj} is satisfied if all $A_j$ are equal.

\item If condition \eqref{hyp:conj} is not true there may be odd degree centralizers.  For example if $A_1=0$ then  $p_1^t$ is a centralizer of $H_0$ for all $t$.

\item Even if condition \eqref{hyp:conj} is satisfied there may be exotic, indecomposable generators of $P$ if the $A_j$ are not the same.  For example in the case $n=2$ if $H_0$ is  the Henon-Heiles Hamiltonian with $A=1$, $B=9$ then $\zbar_1^3 z_2$ is an indecomposable  centralizer.

 \end{enumerate}
}
 \end{rmk}

See \cite{bcu} section 1.5 for related results.   For simplicity we prove Theorem \ref{thm:conj} for the $4$ dimensional Henon-Heiles Hamiltonian.  The proof of the  general case is identical.

 \nd We use Hamilton's equations (see \eqref{def:hameq}, \eqref{Newton}):
 
 \begin{equation}
 \label{H0eqofmotion}
 \dot{q}_j = p_j, \quad j=1,2, \quad 
  \dot{p}_1 = - Aq_1, \quad \dot{p}_2 = - B q_2.
 \end{equation}

 \nd In  terms of the complex variables \eqref{cplxHamgenerators} Hamilton's equations become:

\begin{equation}\label{hameqforEigen}
\dot{z}_1 = -i\sqrt{A} z_1, \quad \dot{z}_2 = -i\sqrt{B} z_2 \end{equation}

$$ \dot{\zbar}_1 = i \sqrt{A}\zbar_1, \quad \dot{\zbar}_2 = i \sqrt{B}\zbar_2. $$

 It follows that the solution to Hamilton's equation, \eqref{def:hameq} is the flow

 $$\varPhi_t(z,\zbar) = ( e^{-i\sqrt{A}t}z_1,e^{-i\sqrt{B}t}z_2, e^{i\sqrt{A}t}\zbar_1, e^{i\sqrt{B}t}\zbar_2).$$

 A monomial $m=z_1^{I_1}z_2^{I_2}\zbar_1^{J_1}\zbar_2^{J_2}$ evaluated on the flow has the form

 $$e^{i (k \sqrt{A} + \ell \sqrt{B})t} z_1^{I_1}z_2^{I_2}\zbar_1^{J_1}\zbar_2^{J_2}=e^{i[m]t}z_1^{I_1}z_2^{I_2}\zbar_1^{J_1}\zbar_2^{J_2},$$

  \nd  where $k= J_1-I_1, \ell = J_2-I_2$.

For monomials $m_1,  m_2$ we say $m_1 \sim m_2$ if $[m_1]=[m_2]$, where $[m]$ is defined in Definition \ref{def:[m]}.

 With these preliminaries we can now prove Theorem \ref{thm:conj}.

{\em Proof of Theorem} \ref{thm:conj}.  Let $K$ be a centralizer of $H_0$.  Then $s(t)=K(\varPhi_t(z,\zbar))$ must be constant, \eqref{chainrule}.  Since $K$ is a finite homogeneous polynomial we may write $K$ in terms of equivalence classes of monomials.  Specifically 

$$K=\sum_{j=1}^{\ell} K_{[m_j]},$$ where $K_{[m_j]}$ is the sum of all the monomials,  $z_1^{I_1}z_2^{I_2}\zbar_1^{J_1}\zbar_2^{J_2}$ in $K$ with $[z_1^{I_1}z_2^{I_2}\zbar_1^{J_1}\zbar_2^{J_2}] = [m_j]$ 
and $\ell$ is the number of equivalence classes of monomials in $K$.

If $K$ is a centralizer of $H_0$  then $$s(t) = \sum_{j=1}^{\ell}  e^{i  [m_j]t}K_{[m_j]}
$$will be constant.

 We show that the terms with non-constant exponential are zero if $s(t)$ is constant. With out loss of generality we may assume $[m_j] \neq 0$.

  Since $s(t)$ is a constant of motion,  $\dfrac{d^r s}{dt^r}=0 $ for all $r>0$.

  Write out $\ell$ derivatives.  We obtain $\ell$ equations:

  $$
  \dot{s}(t) = (i[m_1]) e^{i[m_1]t}K_{[m_1]}+ (i[m_2])e^{i[m_2]t}K_{[m_2]} + \cdots  =0$$
  $$  \ddot{s}(t) = (i[m_1])^2 e^{i[m_1]t}K_{[m_1]}+ (i[m_2])^2e^{i[m_2]t}K_{[m_2]} + \cdots  =0 $$
 $$  \vdots $$
 
  The system of equations may be written as 
  $$A\cdot \nu=0, \quad \nu=( e^{i  [m_1]t}K_{[m_1]},  e^{i  [m_2]t}K_{[m_2]}, \cdots)$$

  \nd where $A$ is a Vandermonde matrix with nonzero determinant.  This implies the vector $\nu$ must be zero. This proves {\it (ii) }.
  
  \nd For $m= z_1^{I_1}z_2^{I_2}\zbar_1^{J_1}\zbar_2^{J_2}$ $s(t) =e^{i[m]t}m $.  By \eqref{chainrule} $ad(H_0)(K) = s^{\prime}(0) = i[m]m$ from which {\it (i)} and {\it (iii)} follow.

\nd  Part {\it (iv)} is a consequence of Proposition \ref{prop:derivation}.  Part {\it (v)} is immediate since $\sum k_i $ equals the degree of the monomial mod $2 $.  Part {\it (vi)} follows since $[m]=0$ if and only if the number of conjugate $z's$ is the same as the number of non conjugate $z$'s. By Corollary \ref{cor:5.10} such monomials are monomials  on the generators of $\text{ker}(ad(H_0)\Big|_{\ch_0}$.  The assertion that $E_1^{2p+1,1}=0$ follows since $[m]=0 \Rightarrow m=0$ if $m \in \ch_{2p+1}$.

  \qed

\begin{cor}{\rm
	
	The basis of the $\ker\{ad(H_0):\ch_0 \to \ch_0\}$ ( $=E_1^{0,0}$ ) is given by
	$$\{ z_i\overline{z}_j | A_i=A_j\}.$$
}
\end{cor}

 \pf  These are the quadratic monomials with $[m]=0$.
 
 \qed

 \begin{cor}\label{cor:5.10}	{\rm

 For $n=2$	the real centralizers of $H_0$ are freely generated.  Specifically

 	\begin{enumerate} 
 		
 		\item If	$A=B$

 		\begin{enumerate}
 			\item $( p_i^2 + Aq_i^2) = z_i \zbar_i, \quad i=1,2$

 			\item $ ( p_1p_2 +A q_1q_2) = \frac{1}{2}( z_1 \zbar_2 +\zbar_1 z_2  )$

 			\item $(p_1q_2 - p_2q_1)= \frac{1}{2 i \sqrt{A}}( z_1 \zbar_2- \zbar_1z_2  )$
 		\end{enumerate}

 		\item If $A \neq B$

 		\begin{enumerate}
 			\item $p_2^2 + Bq_2^2 = z_2\zbar_2$

 			\item $ p_1^2 + Aq_1^2 = z_1\zbar_1$

 		\end{enumerate}

 	\end{enumerate}
}
 	
 \end{cor}

 The previous results apply more generally.  A grading is defined on a polynomial algebra, $\mathcal{L}$ in appendix \ref{sec:hamiltonian}.  Filtration $i-2$, denoted $\mathcal{L}_{i-2}$ is the set of homogeneous polynomials of degree $i$.   In particular $\mathcal{L}_{-2}$ are the constants and $\mathcal{L}_{-1}$ are the linear polynomials.  Assume $\mathcal{L}_{-1}$ has basis $\{x_1, \cdots x_{2n}\}$.  Let $\delta:\mathcal{L}_{-1} \to \mathcal{L}_{-1}$ be a linear map (not necessarily induced by an adjoint representation).  We extend $\delta$ to $P(\mathcal{L}_{-1})$,  the polynomial algebra generated by $\mathcal{L}_{-1}$, by the Leibniz rule:$$\delta(xy) = \delta(x)y+x\delta(y).$$

\nd  There are the eigenvectors,   $\{ v_1, \cdots v_{2n}\}$ of $\delta\Big|_M$
   with eigenvalues $\{e_1, \cdots, e_{2n}\}$  \begin{equation}\label{darboux}\delta(v_i)=e_i v_i \mbox{ for } i= 1, \cdots ,2n.\end{equation}

\begin{assum} \label{assumtion}{\rm
	Assume the vectors $\{v_1, \cdots, v_{2n}\}$ form a basis of $\mathcal{L}_{-1}$.  
}
\end{assum}

\nd Generators  $\{v_1, \cdots, v_{2n}\}$ of $\mathcal{L}_{-1}$ which satisfy \eqref{darboux} are said to satisfy the Darboux condition.

\nd Assuming \eqref{assumtion} the polynomial algebra $P(\mathcal{L})$ has a basis of eigenvectors.   Specifically the monomials 
$m=v_1^{I_1}v_2^{I_2}\cdots v_{2n}^{I_{2n}}$ freely generate $P(\mathcal{L})$.   Since $\delta$ is a derivation we have 
\begin{equation}\label{deltam}\delta(m) = [m] m \mbox{ where } [m]=\sum_{j=1}^{2n} e_j I_j.\end{equation}

  Notice that $[m]$ is logarithmic, i.e. $[mn]=[m]+[n]$ for $m,n$ monomials in $v_i$.

\begin{exam} {\rm

	\nd An example  is given by a Hamiltonian vector field  on $2n$ canonical variables, $\mathbf{x}=\{x_1, \cdots x_{2n}\}=\{q_1, \cdots q_n, p_1, \cdots p_n\}$ (see  \ref{eqofmotion}) with  $\delta=  ad(H_0)$ for $H_0$ a quadratic Hamiltonian.
	Specifically if $H_0$ is the quadratic part of the Henon-Heiles Hamiltonian then
	the vectors $z_1,\overline{z}_1, z_2,\overline{z}_2$ of  \eqref{cplxHamgenerators}, are eigenvectors of $ad(H_0)$  with eigenvalues $$\{-i \sqrt{A}, i \sqrt{A},-i \sqrt{B}, i \sqrt{B}\}.$$ respectively.  These eigenvectors satisfy Assumption \ref{assumtion} if $AB \neq 0$.
}
	\end{exam}

The monomials with zero eigenvalue (e.g. monomials, $m$ with $[m]=0$) are constants of $\delta$.  In order to prove that these monomials generate the constants we have to specialize to the case of a Hamiltonian vector field.  For this case we show that the basis of eigenvectors provides a convenient formulation of Hamilton's equations.

\begin{prop}\label{prop:genham}{\rm If the eigenvectors of $ad(H_0)$ satisfy assumption \ref{assumtion} then the Hamilton differential equation has the form
	$$\dot{v}_i = e_i v_i, \quad i=1, \cdots, 2n.$$
}
	
\end{prop}

\pf  The eigenvector, $v_i$ is a polynomial in $\mathbf{x} = (\mathbf{q},\mathbf{p})$:
$$v_i=v_i(\mathbf{x}).$$   Then
$$\dot{v}_i= \dfrac{d}{dt}(v_i(\mathbf{x}(t)) =\dfrac{d v_i}{d\mathbf{x}} \cdot \dfrac{d\mathbf{x}}{dt} = \langle (v_i)_{\mathbf{x}}, JH_{\mathbf{x}} \rangle = \{ H,v_i\} = e_i v_i.$$
 
 \qed

\nd The equations \eqref{hameqforEigen} are the  special case of  Proposition \ref{prop:genham} for the quadratic term of the Henon-Heiles equation.

\nd Proposition \ref{prop:genham}  allows a description of the solutions to Hamilton's equation and a generalization of Theorem \ref{thm:conj}. 

\nd In light of Proposition \ref{prop:genham} the solution to Hamilton's equation has flow given by

$$\Phi_t(m) = e^{[m] t} m.$$

\nd The argument to prove Theorem \ref{thm:conj} now applies to a general quadratic Hamiltonian satisfying \eqref{assumtion}. In principle this enables the calculation of the $E_1$ term of the spectral sequence of a Hamiltonian with general quadratic term satisfying assumption \eqref{assumtion}.

\begin{exm} {\rm 
	
	We determine the  centralizers of 
	$$H=(p_1q_2-q_1p_2) + A(p_1^2-q_2^2).$$.
	
	\begin{enumerate}
		
		\item  $$ H_x= \Big(-{\it p_2},\quad 2\,A{\it q_2}+{\it p_1}
		\, ,\quad 2\,A{\it p_1}+{\it q_2}
		\, ,\quad -{\it q_1}\,
	\Big)
		$$

		\item $$J  \cdot H_x=\Big(2\,A{\it p_1}+{\it q_2}\,, \quad -{\it q_1}
		\, , \quad {\it p_2}\, ,\quad -(2\,A{\it q_2}+{\it p_1})\, \Big).$$

		\item With respect to the basis $\mathbf{x}=\{q_1,q_2,p_1,p_2\}$  the matrix of the map $ad(H):\mathcal{L}_{-1}\to \mathcal{L}_{-1}$ is:
		
		$$M= \left[ \begin {array}{cccc} 0&-1&0&0\\ \noalign{\medskip}1&0&0&-2\,A
		\\ \noalign{\medskip}2\,A&0&0&-1\\ \noalign{\medskip}0&0&1&0
		\end {array} \right]. $$
		
		\vspace{5 \jot}
		
	\nd 	The eigenvectors of $M$ are the columns of the following array:
		
		$$\left[ \begin {array}{cccc} -1&-1&1&1\\ \noalign{\medskip}-{(2\,A+1){\frac {1}{\sqrt {-2\,A-
						1}}}}&{(2\,A+1){\frac {1}{\sqrt {-2\,A-1}}}}&-\sqrt {2\,A-1}&\sqrt {2
			\,A-1}\\ \noalign{\medskip}-{(2\,A+1){\frac {1}{\sqrt {-2\,A-1}}}}&{(2
			\,A+1){\frac {1}{\sqrt {-2\,A-1}}}}&\sqrt {2\,A-1}&-\sqrt {2\,A-1}
		\\ \noalign{\medskip}1&1&1&1\end {array} \right] .$$

	\nd 	Specifically the eigenvectors are:

		\begin{enumerate}
			
			\item $$v_1= -q_1 -i \sqrt{2A+1} q_2-i \sqrt{2A+1} p_1+p_2, \mbox{ eigenvalue } i \sqrt{2A+1}$$

			\item $$v_2 =  -q_1 +i \sqrt{2A+1} q_2+i \sqrt{2A+1} p_1+p_2, \mbox{ eigenvalue } -i \sqrt{2A+1}$$

			\item $$v_3= q_1-\sqrt{2A-1}q_2  +\sqrt{2A-1} p_1 +p_2, \mbox{ eigenvalue} \sqrt{2A-1}$$

			\item $$v_4= q_1+\sqrt{2A-1}q_2  -\sqrt{2A-1} p_1 +p_2, \mbox{ eigenvalue} -\sqrt{2A-1}$$
			
		\end{enumerate}

		\item The determinant of the matrix $[v_1,v_2,v_3,v_4]$ is 
		
		$$16\,\sqrt {-2\,A-1}\sqrt {2\,A-1}.$$

		So the vectors $\{v_1, \cdots , v_4\}$ satisfy assumption \eqref{assumtion} if $A\neq \pm \frac{1}{2}$.

		\item For a monomial $m:= v_1^{I_1}v_2^{I_2}v_3^{I_3}v_4^{I_4}$ the bracket function is given by

		$$[m]=i I_1 \sqrt{2A+1} -i I_2\sqrt{2A+1} +I_3 \sqrt{2A-1} -I_4 \sqrt{2A-1}$$

		and 
		$$ad(H)(m) = [m]m.$$

	\end{enumerate}

}
	\end{exm}

 Finally we note that these considerations extend to the field or rational functions, $Q = \{ \dfrac{f}{g} | f,g \in P\}$.  The formula for $\delta$ on monomials in $v_i$ extends to negative exponents.  Specifically if $m$ and $n$ are monomials in eigenvectors
$$\delta(\dfrac{m}{n}) = \dfrac{m}{n} ([m]-[n]).$$

\nd In particular the ratio of two eigenvectors with the same eigenvalues are constants for the differential field, $Q$ with derivation $\delta$.

In the following examples we use the spectral sequence to find centralizers of the Henon-Heiles Hamiltonian  $H$ (Definition \eqref{def:henon-heiles}). We find centralizers  of $H$ which have non-trivial quadratic term and are of degree less than or equal to $4$.  In Corollary \ref{cor:5.10}  we  determined $E_1^{0,0}= \text{ker} (ad(H_0)\Big|_{\ch_0}$.  We make use of MAPLE to study the problem of extending the generators of $E_1^{0,0}$ given in Corollary \ref{cor:5.10}, $(i)$ and $(ii)$ case $(a)$.

 \begin{expls}\label{exm:314}{\rm 
 	
 	  In this example we give a detailed analysis of how the spectral sequence may be used to study the problem of extending the class 
 	
 	$$K_0= \frac{3}{2}(4B-A)( Bq_2^2+p_2^2)$$

 	\nd to a centralizer of $H$   (the coefficient appears to minimize messy coefficients that appear in denominators).  We regard $A$ and $B$  to be indeterminate and determine which $L$ allow centralizers which are independent from $H$. In particular since $A$ and $B$ are independent, condition \eqref{hyp:conj} is satisfied.   The reader is referred to Appendix \ref{app:ss} for guidance and to \cite{bc} examples 4.1.1 and 4.1.2.

 We first compute $d_1:E_1^{0,0} \to E_1^{1,1}$
 We know from Theorem \ref{thm:conj} that  $E_1^{1,0}=0$ but we have to keep 
 	track of how classes are killed in $E_1^{1,0}$ in order to compute $d_2$.  Specifically we compute

 	$$d_1(K_0) = \{H,K_0\} = \{H_0,K_0\} + \{H_1,K_0\} = 6\,L{\it p_2}\,{\it q_1}\,{\it q_2}\, \left( A-4\,B \right) .$$

 	Since $E_1^{1,0}=0$ there must be $K_1\in \ch_1$ such that $$ad(H_0)(-K_1)=6\,L{\it p_2}\,{\it q_1}\,{\it q_2}\, \left( A-4\,B \right).$$
 	Theorem \ref{thm:conj} (i) tells us how to find $K_1$.  Namely if $6\,L{\it p_2}\,{\it q_1}\,{\it q_2}\, \left( A-4\,B \right)$ is written in terms of monomials $m_i$ in the complex variables $z_i ,\zbar_i$ then $ad(H_0)^{-1}(m_i) = \dfrac{m_i}{[m_i] i}$.  With MAPLE's help one finds
 
 $$-K_1=6\, \left(  \left( B{{\it q_2}}^{2}-{{\it p_2}}^{2} \right) {\it q_1}+{
 	\it p_2}\,{\it q_2}\,{\it p_1} \right) L.$$

\nd So $K=K_0+K_1$ has the property that $ad(H)(K) \in \ch_2$.  In the language of spectral sequences, the class $[K_0] \in E_1^{0,0}$ survives to  $[K_0+K_1]  \in E_2^{0,0}$.

\nd  We now compute $d_2: E_2^{0,0} \to E_2^{2,1}$.  As before we compute $d_2$ by evaluating

 $$d_2([K_0])=ad(H)(K)=\overbrace{\{H_0.K_0\}}^{=0}+\overbrace{\{H_1,K_0\}+\{H_0,K_1\}}^{=0}+\{H_1,K_1\} \in \ch_2.$$

\nd  Referring to \eqref{nzdefs} we have

$$E_2^{2,1} = \dfrac{F^2 \ch}{ad(H)(Z_1^1) + F^3 \ch}.$$

\nd where $Z_1^1=\{m=m_1+m_2| ad(H)(m) \in F^2 \ch\}$.   For $m \in Z_1^1$ we have
$$\{H_0+H_1,m_1+m_2\} = \{H_0,m_1\} +\{H_1,m_1\} + \{H_0,m_2\}.$$

\nd The term $\{H_0,m_1\}$ must be zero since $ad(H)(m) \in F^2$.  From Theorem \ref{thm:conj}, \eqref{hyp:conj} we have  $m_1=0$.  So we conclude that $E_2^{2,1} = \dfrac{ F^2\ch}{\{H_0,\ch_2\}+F^3 \ch}$.

\nd Using Theorem \ref{thm:conj} one  readily checks (using MAPLE of course) that $d_2([K_0])=\{H_1,K_1\}=0$ in $E_2^{2,-1}$.  Namely, as we did for $d_1$ we expand $\{H_1,K_1\}$ in complex coordinates and use Theorem \ref{thm:conj} (i) to find a  $K_2$
  such that  $\{H_0,-K_2\}= \{H_1,K_1\}$.  Then  $$K^{\prime}=K_0+K_1+K_2$$ has the property that $ad(H)(K^{\prime})\in \ch_3$.

\nd   We will not burden  the reader with the formula for $K_2$ other than to observe that it is quite complicated.  The next step is to compute $d_3([K_0]) = \{H_1, K_2\}$.

\nd  MAPLE does provide a $K_3$ such that $ad(H_0)(K_3) = \{H_1,K_2\}$ which can be used to perturb $K^{\prime}$ to a polynomial that is a centralizer of $H$ modulo filtration $4$.  However for some values of $L$ we can actually perturb $K^{\prime}$ to $K^{\prime \prime}$ so that $d_3([K_0]) =0 \in \ch_3$.  The resulting polynomial,  $K^{\prime \prime}$ is a centralizer of the Henon-Heiles Hamiltonian.

\nd  To this end we recall that $d_3([K_0]=ad(H)(K^{\prime})=\{H_1,K_2\}$ takes values in $$E_3^{3,1}=\dfrac{F^3\ch}{ad(H)(Z_2^1)+F^4\ch},$$

  \nd  where $m=m_1+m_2+m_3 \in  Z_2^1$ means
\begin{multline}\label{adhm}  ad(H)(m)=\{H_0+H_1,m\}= \\\{H_0,m_1\}+\{H_1,m_1\}+\\
\hspace{-2 in}\mbox{ (i) }\hspace{1.6 in}\{H_0,m_2\}+\{H_1,m_2\}+ \\
\hspace{-2.7in} \mbox{(ii)} \hspace{2 in} \{H_0,m_3\}\\
\end{multline}

\nd is in $F^3 \ch$.  The implications of $ad(H)(m) \in F^3 \ch$  are:

\begin{itemize}
\item  $\{H_0,m_1\}+\{H_1,m_1\}=0 $

\item  $\{H_0,m_2\}=0$,  which implies   $m_2 \in E_1^{2,0}$ and $\{H_1,m_2\}=d_1(m_2)$.

\item $\{H_0,m_3\} = d_0(m_3)$.

\end{itemize}

    As mentioned above,  $d_3([K_0])$ is zero in $E_3^{3,1}$ if there is an $m_3$ (which we called $K_3$) such that  $\{H_1,K_2\}=\{H_0,m_3\}$.  But 
  $d_3([K_0])$ is also zero if there is $m_2 \in E_1^{2,-2}$ such that $\{H_1,m_2\} =\{H_1,K_2\}$.  From Theorem \ref{thm:conj} we know that $E_1^{2,-2}$ is generated by $(z_1\overline{z}_1)^2, (z_2\overline{z}_2)^2$ and $z_1\overline{z}_1 z_2\overline{z}_2$.  Modulo $q_1$ and $q_2$ these are $p_2^4,\quad  p_1^4$ and $p_1^2p_2^2$ respectively. In $K_2$  we notice that the monomial $p_2^4$ occurs with coefficient $\dfrac{9 L^2}{16 B^2}$, the monomial $p_1^4$ has coefficient $0$ and $p_1^2p_2^2$ has coefficient $\dfrac{3}{4}\,{\frac {L \left( 4\,AL+2\,BL-A \right) }{BA \left( A-B \right) }}$.

\nd Since the image of $ad(H_1)$ is zero modulo $q_1$ and $q_2$ this discussion motivates perturbing $K_2$ by subtracting $m_2$ where

  \nd    \begin{equation}\label{someterms} m_2 =-\dfrac{9 L^2}{16 B^2} (z_2\overline{z}_2)^2 + \dfrac{3}{4}\,{\frac {L \left( 4\,AL+2\,BL-A \right) }{BA \left( A-B \right) }}z_1\overline{z}_1z_2\overline{z}_2. \end{equation}

 \nd    Here is $K_2-m_2$:

 $-6\,{\frac {L \left( -1/2\,L \left( A-B \right) {{\it q_2}}^{4}+
 		\left(  \left( A-B \right)  \left( L-1/2 \right) {{\it q_1}}^{2}+3\,
 		\left( L-1/6 \right) {{\it p_1}}^{2} \right) {{\it q_2}}^{2}-6\,{\it p_2
 		}\, \left( L-1/6 \right) {\it p_1}\,{\it q_1}\,{\it q_2}+3\,{{\it p_2}}^{2
 		} \left( L-1/6 \right) {{\it q_1}}^{2} \right) }{2\,A-2\,B}}$

\nd Quite a mess.  But the differential, $d_3([K_0])\in E_3^{3,1}$ is more tractable.  Letting $K^{\prime\prime} = K^{\prime} -m_2$ we have

 \begin{equation}\label{Bq_2example}
 \{H, K^{\prime\prime}\}=-3\,{\frac {L \left( {\it p_1}\,{\it q_2}-{\it p_2}\,{\it q_1} \right) 
 		\left( 6\,L-1 \right)  \left( 2\,L{{\it q_1}}^{2}-L{{\it q_2}}^{2}-{{
 				\it q_1}}^{2} \right) {\it q_2}}{A-B}} \in \ch_3.
 			\end{equation}

 \nd The process stops with $4-$th degree centralizer if \eqref{Bq_2example} is zero.   This is the case if  $L= \dfrac{1}{6}$ or $0$.  If $L=\dfrac{1}{6}$ we have (compare with \cite{cf} Example 1).

 $$K= \dfrac{3}{2}(4B-A)(Bq_2^2+p_2^2) +Bq_1q_2^2+p_2(q_2 p_1 - q_1p_2) + \dfrac{1}{6}q_2^2(q_1^2 + \dfrac{1}{4}q_2^2).$$

 If $L=0$ we obtain the centralizer $$p_2^2+Bq_2^2.$$

\nd By similar methods one can compute extensions of $$Aq_1^2+p_1^2$$ to a centralizer.  One finds that  the analogue of formula \eqref{Bq_2example} is zero only if $L=0.$   The resulting centralizer is $$ Aq_1^2+p_1^2 + \dfrac{2}{3} q_1^3.$$ 

The details are left to the reader.

We conclude that without any restrictions on $A$ or $B$ all centralizers of degree $4$ or less are generated by:

\begin{enumerate}[label=(\roman*)] \label{enu:AnotB}

 \item $L=0$

 $$Aq_1^2+p_1^2 + \dfrac{2}{3} q_1^3 \mbox{ and } Bq_2^2+p_2^2.$$

 \item $L =1/6$

   $$\dfrac{3}{2}(4B-A)(Bq_2^2+p_2^2) +Bq_1q_2^2+p_2(q_2 p_1 - q_1p_2) + \dfrac{1}{6}q_2^2(q_1^2 + \dfrac{1}{4}q_2^2).$$
\end{enumerate}
}
\end{expls}

\begin{expls}{\rm

If $A=B$ there are two more generators in $E_1^{0.0}$ (see Corollary  \ref{cor:5.10}).

$$ p_1p_2 +A q_1q_2  \mbox{ and }  p_1q_2 - p_2q_1.$$

We search for centralizers of $H$ by computing the differentials on an arbitrary linear combination of generators of $E_1^{0,0} $ written in complex coordinates.\\

$K_0=az_1\zbar_1 + b z_2 \zbar_2 + c \zbar_1 z_2 +d z_1 \zbar_2 =$
$$-i \left( {\it p_1}\,{\it q_2}-{\it p_2}\,{\it q_1} \right)  \left( c-d
\right) \sqrt {B}+a{{\it p_1}}^{2}+{\it p_2}\, \left( c+d \right) {\it 
	p_1}+b{{\it p_2}}^{2}+$$$$B \left( a{{\it q_1}}^{2}+{\it q_2}\, \left( c+d
\right) {\it q_1}+b{{\it q_2}}^{2} \right. $$

  We assume both $c$ and $d$ are not zero (which is the case of the previous example).  Proceeding as in Example \ref{exm:314} we compute  $d_1= \{H,K_0\} \in \ch_1$ and find it is not zero.  Therefore there is a unique class, $K_1$ such that $d_2=\{H,K_0+K_1\} \in \ch_2$.  But in contrast  to the previous example   this class is not zero in $E_2^{2,1}$.  It is a sum of $6$ non zero classes, a typical one is given below:

$$ \left( -\dfrac{i}{12}Ld+\dfrac{i}{2}{L}^{2}d+\dfrac{i}{12}{L}^{2}c-\dfrac{i}{2}Lc \right) {\it \zbar_1}\,{{
		\it z_2}}^{2}{\it \zbar_2}.$$

	Notice that this class is indeed not zero in $E_2^{2.1}$ since $[{\it \zbar_1}\,{{
			\it z_2}}^{2}{\it \zbar_2}]=0$ ($[{\it \zbar_1}\,{{
			\it z_2}}^{2}{\it \zbar_2}]$ is defined in Definition \ref{def:[m]}).

	These obstructions are killed by setting $a=b$ and $c=d$ and $L=1$.   One finds that not only is $d_2$ now zero, but that $\{H,K_0+K_1\}$ is zero in $\ch$.  And we have the following family of centralizers of $H$.

\begin{equation}\label{K0K1}K_0+K_1= a\left(  \left( {{\it q_1}}^{2}+{{\it q_2}}^{2} \right) B+{{\it p_1}}^{2}
+{{\it p_2}}^{2} \right) +2\,d \left( B{\it q_1}\,{\it q_2}+{\it p_1}\,{
	\it p_2} \right) 
	+ \end{equation}$$2\,d{{\it q_1}}^{2}{\it q_2}+2\,a{\it q_1}\,{{\it q_2}}^{2}+\dfrac{2}{3}\,d{{\it q_2
	}}^{3}+\dfrac{2}{3}\,a{{\it q_1}}^{3}.$$
	
	This completes the description of all the centralizers of the Henon-Heiles Hamiltonian with quadratic lowest term of degrees bounded by $4$.   Letting $a=\dfrac{1}{4}, \quad d=-\dfrac{1}{4}$ in \eqref{K0K1} agrees with  example 4.1.3 in \cite{bc}.

}

\end{expls}

 \appendix

 \section{Spectral Sequences}\label{app:ss}

 In this appendix we discuss  spectral sequences that arise from  filtered cochain complexes.  We limit our attention to the cochain complex of  a filtered graded linear map.  An interested reader is directed to  \cite{keymac} for more details.  Throughout   $\,\mbk\,$ denotes a commutative ring with
 multiplicative identity $\,1 \neq 0\,$ and all modules, $M$ are assumed
 (left) $\, \mbk$-modules.

 A spectral sequence is  a collection  $\{E_r\}$ of $\, \mbk $ modules equipped with maps $d_r:E_r \to E_r$ with $d_r^2=0$.  Such a module is called a {\it differential module}.   The {\it homology} of a differential module is the $\, \mbk  $ module $H(E_r,d_r)=\text{ker}(d_r) /\text{Im}(d_r)$.  The sequence of differential modules,  $\{E_r,d_r\}.$ is called a spectral sequence if there is an isomorphism 
 $H(E_r,d_r) = E_{r+1}$.   The terms in the spectral sequence are often viewed as the pages of a book and  $E_r$ is often called the $r-$th page of the spectral sequence.

 There is a  spectral sequence induced by  a filtered cochain complex.  We describe the particularly simple filtered cochain complex  that arises from a graded linear map.

 Let  $\,L\,$  be an  $\,\mbk$  module.  A family $\,\{F_pL\}_{p \geq 0}\,$ of submodules is an {\it increasing filtration}\, of $\,L\,$ if 
 \begin{equation} \label{infil}
 	i < j  \qquad \Rightarrow \qquad F_{i}L \subset F_{j}L;
 \end{equation}
 a family $\,\{F^{p}L\}_{p \geq 0}\,$  is a {\it decreasing filtration}\, of $\,L\,$ if
 \begin{equation} \label{defil}
 	i < j \qquad \Rightarrow \qquad F^{i}L \supset F^{j}L.
 \end{equation} 
 The {\it associated graded}\, ({\it module} of the family)\, in the increasing case refers to the collection of quotients 
 \begin{equation} \label{agone}
 	\{{\rm Gr}_{p}L\}_{p \geq 0}, 
 \end{equation}
 where
 \begin{equation} \label{agtwo}
 	{\rm Gr}_{p}L := F_{p}L/F_{p-1}L \qquad \text{when} \qquad p > 0
 \end{equation}
 and
 \begin{equation}  \label{grz}
 	{\rm Gr}_{0}L := F_{0}L.
 \end{equation}
 The {\it associated graded}\, ({\it module} of the family)\, in the decreasing case refers to the collection of quotients 
 \begin{equation} \label{five}
 	\{{\rm Gr}^pL\}_{p \geq 0}, 
 \end{equation}
 where
 \begin{equation} \label{agsix}
 	{\rm Gr}^{p}L := F^{p}L/F^{p+1}L.
 \end{equation}
 
 We obtain filtrations  when $\,L\,$   {\it graded,}\, that is to say when $\,L\,$ is a module of the form
 \begin{equation}  \label{asgrinc}
 	L := \ts \bigoplus_{j \geq 0}L_j.
 \end{equation}
 the {\it associated} ({\it increasing}) {\it filtration}  $\,\{F_pL\}_{p \geq 0}\,$ of $\,L\,$ is that defined by
 \begin{equation} \label{asgrincr}
 	F_{p}L := \ts \bigoplus_{0 \leq j \leq {p}}L_{j},
 \end{equation}
 and the {\it associated} ({\it decreasing}) {\it filtration}   $\,\{F^{\bf p}L\}_{p \geq 0}\,$ is that defined by
 \begin{equation}\label{asgrdcr}     F^{p}L := \ts \bigoplus_{j \geq p}L_{j}.
 \end{equation}

 Elements of any particular $\,L_{j}\,$  will be said to  be at (or to have) {\it grading level} $\,{ j}$. Elements of $\,F_pL\,$ or $\,F^pL\,$ will be said to be at {\it filtration level} $\,p$. 
 When the grading/filtration alternative should be clear  from context we will simply refer to elements ``at level $j$.'' 
 
 For our purposes the advantage in assuming\, (\ref{asgrinc})\, arises from the identifications
 \begin{equation} \label{asgrid}
 	{\rm Gr}_{p}L \simeq L_p \simeq {\rm Gr}^{p}L.
 \end{equation}
 In particular,  in the increasing case any  $\,\ell \in F_pL\,$ can be written uniquely as a finite sum
 \begin{equation} \label{asgrfour}
 	\ell = \ell_{0} + \cdots + \ell_{r} + \cdots + \ell_{p}, \qquad \text{where} \qquad \ell_{ j} \in L_{ j}. \ \  
 \end{equation}
 and the subscripts are strictly increasing.   The isomorphism \eqref{asgrid} associates the equivalence class $\, [\ell] \in {\rm Gr}_pL$ with $\, \ell_p \in L_p$.   The analogue  in the decreasing case  of a class in $F^pL$ is a finite sum
 \begin{equation} \label{asgrfive}
 	\ell = \ell_{p} + \cdots + \ell_{ r} + \cdots + \ell_{s}, \qquad \text{where} \qquad \ell_{ j} \in L_{ j} .\ \ \ 
 \end{equation}
 and the subscripts are again strictly increasing.  In both cases we refer to  $\,\ell_{p}\,$ as the {\it leading term}\, of $\,\ell$, which we often write as $\,\overline{\ell}$.

 \begin{exams} \label{Lexams}
 	{\rm The following examples of graded modules played key roles in Sections \ref{section:vectorfields} and \ref{sec:appHam} .  
 		\begin{itemize}
 			\ia Assume the total degree ordering\footnote{a/k/a the ``degree lex(icographical) ordering"} on the collection of monomials in the variables $\,x,y$,\, in which case
 			\begin{equation*} 
 				\hspace{-.4in}\text{(i)} \hspace{1in}
 				1 \prec x \prec y \prec x^2 \prec xy \prec y^2 \prec x^3 \prec x^2y \prec \cdots \ .
 			\end{equation*}   Take $\,L = \oplus_{p \geq 0}L_p$,\, where $\,L_p\,$ is the free $\,A$-module generated by the $\,p^{th}$-term $\,x^ry^s\,$ on the list,\, i.e.
 			\begin{equation*}\hspace{-1.18in}\text{(ii)} \hspace{1.56in}
 				p := \ds  \frac{(r+s)((r+s)+1)}{2} + s. \smallskip
 			\end{equation*}
 			
 			The mapping $\,x^ry^s \mapsto p$,\,  is an order preserving isomorphism between
 			the multiplicative monoid $\,\{x^ry^s\}_{(r,s) \in \nn^2}$,\, ordered as in\, (i),\, and the additive monoid $\,\nn$  $ := \{0,1,2,3,\dots\,\}\,$ with the usual order $\,<$\,. (Note that $\,x^0y^0 = 1 \mapsto 0$.)   If we define  $\alpha(n)=\lfloor\frac{1}{2}(-1+\sqrt{1+8n})\rfloor$ then the inverse is given by $p \mapsto (\alpha(p), \frac{1}{2}(\alpha(p)-1)(\alpha(p))$.
 			
 			This example can be extended, assuming the total degree ordering, to the case of $\,n\,$ variables.   If a different ordering is assumed  the resulting indexing can be inappropriate for our purposes. For example,   for the pure lex ordering
 			\begin{equation*}
 				1 \prec x \prec x^2 \prec \cdots \prec x^n \prec \dots y \prec xy \prec \cdots 
 			\end{equation*}
 			there are infinitely many monomials between $\,x\,$ and $\,y$.  To deal with spectral sequences (as they are treated in this paper) one needs to assume only finitely many predecessors for each monomial.  
 			
 			\ib Fix positive integers $\,m\,$ and $\,n$,\  and  for each $\,p \geq 0\,$ let $\,L_p\,$ denote the collection of $\,m \times 1\,$ column vectors having homogeneous polynomials in $\,\mbk[x_1,x_2,\dots,x_n]\,$  of degree $\,p+1\,$ as entries.  Define $\,L := \oplus_{p \geq 0}L_p$.  This is the grading used in section \ref{section:vectorfields}.
 			\ic Fix a positive integer $\,n\,$ and for each integer $\,p \geq 0\,$ let $\,L_p\,$ denote the collection of homogeneous polynomials in the algebra $\,\mbk[x_1,x_2,\dots,x_{2n}]$ of degree $\,p+2$.  Define $\,L := \oplus_{p \geq 0}L_p$.   This is the grading used in section  \ref{sec:appHam}.
 		\end{itemize}
 	}
 	\end{exams}

 	The assumptions  for the squeal   are as follows: \smallskip
 	\begin{equation}\label{assumptions}\end{equation}
 	\begin{itemize}
 		%\ia $\mcj\,$ is a monoid with an admissible ordering $\,\prec$; \vspfv
 		\ia $M= \oplus_{j \geq 0}M_j\,$ is a graded module with decreasing filtration $\,\{F^{p}M\}_{p \geq 0}\,$   arising as in\, (\ref{asgrincr})  or increasing filtration $\{F_p\}$ as in  \eqref{asgrdcr};\vspfv
 		\ib $N\,$ is a (not necessarily graded) module with decreasing filtration $\,\{F^pN\}_{p \geq 0}$; or increasing filtration $\,\{F_pN\}$ and  \vspfv 
 		\ic $f:M \rightarrow N\,$ is a module homomorphism which preserves the \\  \hspace*{.08in} filtrations, i.e.\
 		satisfies $\,f(F^{p}M) \subset F^{p}N\,$ for all $\,p \geq 0$.\vspfv
 		%\id \ when $\,N\,$ is assumed graded the increasing filtration $\,\{F_{p}N\}_{p \geq 0}\,%$  on\ $N\,$ arises as in\, (\ref{asgrincr}). \medskip
 	\end{itemize}

 		For $\,\mbk$ modules $M$ and $N$ a $\,\mbk$-linear mapping $\,f:M \rightarrow N \,$  can be embedded into the finite cochain complex	i.e. $f$ can be considered as one mapping within the cochain complex

 	\begin{equation}\label{chaincomplex}
 		\mathcal{C}:\cdots \rightarrow 0 \rightarrow 0 \rightarrow C_0 := M
 		\stackrel{f}{\rightarrow} C_1:= N
 		\stackrel{0}{\rightarrow} 0 \rightarrow 0 \rightarrow \cdots\
 		.
 	\end{equation}

 	We now assume $\,M\,$ and $\,N\,$ admit decreasing filtrations satisfying \eqref{assumptions}.  Then $f$ fits into a cochain complex

 	\begin{equation}\label{decreasingfilt}
 		\cdots \rightarrow 0 \rightarrow 0 \rightarrow F^p M 
 		\stackrel{f}{\rightarrow} F^p N
 		\rightarrow 0 \rightarrow 0 \rightarrow \cdots\
 		.
 	\end{equation}
 	(there is a similar filtered cochain complex associated to an increasing filtration). \\

\nd  A filtered cochain complex induces a   spectral sequence:

 	\begin{equation} \label{def:decss}
 		\hspace{-4.5 in} E_r^{\ast,\ast}.\end{equation} in the case of a decreasing filtration and
 	\begin{equation}
 		\label{def:ssinc} \hspace{ -4.5 in}E^r_{\ast,\ast}.\end{equation}
 	in the case of an increasing filtration.\\
 	
 In general the spectral sequence of a filtered cochain complex give information about the homology of the cochain complex.	   In our case the spectral sequences are approximation to:\\

 	\begin{itemize}
 		\item $H^0\mathcal{C} = \text{ker}(f) $ \\
 		
 		\item $H^1\mathcal{C} = \text{coker}(f) $\\
 		\item $0$ otherwise.
 	\end{itemize}

 	The spectral sequence has a natural bigrading, $(p,s)$.  The component  $p$ coming from the filtration and the component $s$  indicates the position in the cochain complex (so in our case  there are only two possibilities). \\
 	
 		We restrict our attention in this appendix  to describing the $E_r$ page of the spectral sequence associated to the short  filtered cochain complexes induced  by a graded linear map  \eqref{decreasingfilt}    (see \cite{keymac} for the general case)\\
 	
 We first consider the case of an increasing  filtration.  	For any non-negative integers $\,p\,$ and $\,r\,$ define \medskip
 	\begin{equation} \label{zdefs}
 		\begin{array}{lllll}
 			Z^r_p & := 
 			& \left\{ \ \ \  
 			\begin{array}{lllllll}
 				\{\,m \in F_{p}M\,:\,fm \in F_{p-r}N\,\} \ \ &\text{if} \ \ 0 \leq r , \vspfv \\
 			
 			\end{array} \right.  \quad
 			\vspfv \\
 			M^r_p & := & \Big\{\ \ \ \ \ \big\{\,m_p \in M_p\,:\, m_p = \overline{m}\, \ \ \text{for some}\ \  \ m \in Z^r_p\big\} .\vspfv \\
 		\end{array}
 	\end{equation}
 	It follows easily from  the filtration-preserving hypothesis that:
 	\begin{equation} \label{sltr}
 		\left\{ \ \ \ \ \ \ \begin{array}{llllll|}
 			
 			{\rm (a)} & Z^{0}_{p} = F_{p}M;\vspfv \\
 			{\rm (b)} & M^{0}_{p} ={\rm Gr}_pM \simeq  M_{p}; \ \ \text{and} \vspfv \\
 			{\rm (c)} & 
 			Z^r_p \subset Z^s_p \ \ \text{if} \ \ r \geq s\,.
 		\end{array} \right.
 	\end{equation}
 	For example, to establish the order reversing correspondence (c)\, simply note that for $\,m \in F_pM\,$ one has
 	\begin{equation*}
 		\begin{array}{lllll}
 			m \in Z^r_p & \Leftrightarrow & fm \in F_{p-r}N  \vspfv \\
 			& \Rightarrow & fm \in F_{p-s}N \ \ \text{(because the filtration is increasing)} \vspfv \\
 			& \Rightarrow & m \in Z^s_p.
 		\end{array}
 	\end{equation*}
 	
 \nd 	Define $\, E^0_{p,0}= {\rm Gr}_pM,  E^0_{p,1}= {\rm Gr}_p N$ and $d_0={\rm Gr}(f)$.  	For $r\geq 1$	define $E^r_{p,0}=M_p^r$ and $E^r_{p,1}=\frac{F_{p}N}{f(Z^{r-1}_{p+r-1}) +  F_{p-1}N},$
 	
\nd  The restriction $\,f\big|_{Z^r_p}\,$ induces a linear mapping
 	\begin{equation} \label{thedrmap}
 		d_r:m_p \in E^r_{p,0} \mapsto [fm] \in E^r_{p-r,1} ,
 	\end{equation}
 	where $\,m \in Z^r_p\,$ is any element satisfying $\,\overline{m} = m_p$,\, and (here and elsewhere) brackets $\,[ \ \, ]\,$ are used to indicate cosets of  the (relevant) quotient.
 	When $\,r > 0\,$ the elements of $\,Z^r_p\,$ are said to ``drop filtration'' (under the mapping $\,f$).
 	
 	%Since definition\, (\ref{thedrmap})\, depends on  both $\,p\,$ and $\,r,\,$ one might prefer replacing the notation $\,d_r\,$ with, say, $\,d^p_r$.  We will not do so: the $\,d_r\,$ notation is standard.
 	
 	\nd It is readily checked that $E^{r+1}_{*,*}= H(E^r_{*,*})$  for $r \geq 0$.

 	The mapping $\,d_r\,$ is easily remembered by means of the commutative diagram
 	\begin{equation} \label{grdiag}
 		\begin{array}{ccc}
 			Z^r_p & \stackrel{f|_{Z^r_p}}{\longrightarrow} & F_{p-r}N \vspfv \\
 			\big\downarrow & & \big\downarrow \vspfv \\
 			M^r_p & \stackrel{d_r}{\longrightarrow} &  
 			\frac{F_{p-r}N}{f(Z^{r-1}_{p-1}) \, + \, F_{p-r-1}N},
 			%\frac{N_{\bf r}}{f(Z^{\bf p}_{{\bf r}^{\rm is}})\cap N_{\bf r}}
 		\end{array} \vspsv 
 	\end{equation}
 	wherein the left vertical mapping is the restriction of the projection of \newline $\,F_{p}M = \ts \!\!\oplus_{0 \leq j \leq {p}}M_{j}\,$ onto $\,M_{p}\,$ and the right vertical mapping 
 	is the usual quotient mapping.  When $\,r = 0\,$  this diagram reduces to
 	\begin{equation} \label{zgrdiag}
 		\begin{array}{rrrcccccccc}
 			& & \ \ \ \ \ \ \ \ F_pM \ \ \    \stackrel{f|_{F_pM}}{\longrightarrow} & F_{p}N \vspsv \\
 			&  & \big\downarrow \ \ \ \ \ \ \ \ \  \, \ \  &   \big\downarrow \vspfv \\
 			{\rm Gr}_pM & = & \ds \frac{F_{p}N}{ F_{p-1}N} \simeq M_p   \ \ \ \    \stackrel{d_0}{\longrightarrow} &  
 			\ds \frac{F_{p}N}{ F_{p-1}N} & = & {\rm Gr}_pN \ \simeq \ N_p\,
 			%\frac{N_{\bf r}}{f(Z^{\bf p}_{{\bf r}^{\rm is}})\cap N_{\bf r}}
 		\end{array} \vspsv 
 	\end{equation}
 	with the understanding that $\,F_{p-1}N := 0\,$ when $\,p = 0\,$.

 	If  $\,M\,$ and $\,N\,$ are graded , and the increasing filtration arises as in\, (\ref{asgrincr}),\,
 	diagram\, (\ref{grdiag})\, can be replaced by
 	\begin{equation} \label{grdiagtwo}
 		\begin{array}{ccc}
 			Z_{\rm p}^{\rm r} & \stackrel{f|_{Z^{r}_{p}}}{\longrightarrow} & F_{p-r}N \vspfv \\
 			\big\downarrow & & \big\downarrow \vspfv \\
 			M^r_{\bf p} & \stackrel{d_{r}}{\longrightarrow} &  
 			\frac{N_{p-r}}{f(Z^{r-1}_{p-1}) \cap N_{p-r}},
 			
 		\end{array} \vspsv 
 	\end{equation}
 	and\, (\ref{zgrdiag})\, by 
 	
 	\begin{equation} \label{grdiagz}
 		\begin{array}{ccc}
 			F_pM & \stackrel{f|_{F_pM}}{\longrightarrow} & F_pN \vspfv \\
 			\big\downarrow & & \big\downarrow \vspfv \\
 			M_p & \stackrel{d_0}{\longrightarrow} &  
 			N_p.
 			
 		\end{array} \vspsv 
 	\end{equation}
 	\nd  It is convenient to display 	the spectral sequence of an increasing filtration as a chart  in the first quadrant of the $(p,q)$ plane. 
 	
 	The $d_2 $ differentials on the $E^2$ page are depicted by SE arrows which drop $2$ filtrations.
 	
 \newpage

 
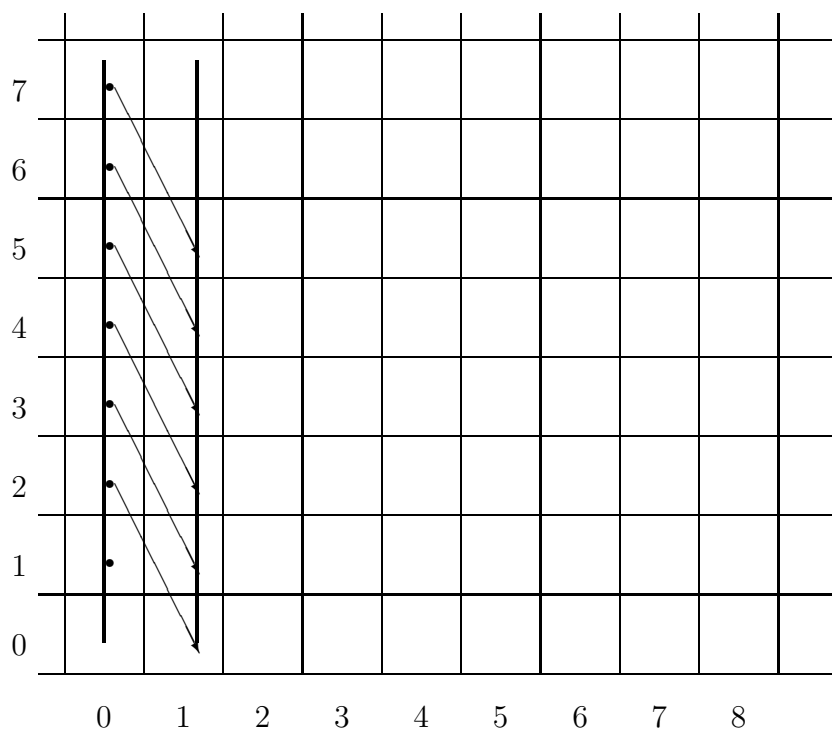
\begin{center}
 	\begin{picture}(420,500)
 		
 		\hspace*{-.6in}
 		\def\mp{\multiput}
 		\def\elt{\circle*{3}}
 		\mp(89,267)(0,30){6}{\vector(1,-2){32}}
 		
 		\mp(60,195)(0,30){9}{\line(1,0){300}}
 		\mp(70,195)(30,0){10}{\line(0,1){250}} \put(142,175){$2$}
 		
 		\put(172,175){$3$} \put(202,175){$4$} \put(232,175){$5$}
 		\put(262,175){$6$} \put(292,175){$7$} \put(322,175){$8$}
 		\put(112,175){$1$} \put(82,175){$0$}
 		
 		\put(50,412){$7$}\put(50,382){$6$}\put(50,352){$5$}
 		\put(50,322){$4$}\put(50,292){$3$}\put(50,262){$2$}
 		\put(50,232){$1$}\put(50,202){$0$}

 		\linethickness{0.35mm}
 		\mp(85,427)(35,0){2}{\line(0,-1){220}}

 		\mp(87,417)(35,0){1}{\elt} \put(120,417)
 		
 		\mp(87,387)(35,0){1}{\elt} \put(120,387)
 		
 		\mp(87,357)(35,0){1}{\elt} 
 		
 		\mp(87,327)(35,0){1}{\elt} 
 		
 		\mp(87,297)(35,0){1}{\elt}

 		\mp(87,267)(35,0){1}{\elt}

 		\mp(87,237)(35,0){1}{\elt}
 		
 		\put(55,160){{\it The } $E^2_{p,q}$ {\it term with differentials from coordinate } $(p,0)$}
 		\put(120,140)	{ {\it to coordinate } $ (p-2,1)$}
 	\end{picture}

 \end{center}

\newpage
  The decreasing case is handled, for the most part, by interchanging sub and superscripts and, to a lesser extent, minus and plus signs, in the development given for the increasing case.  
 %Moreover, in our examples we only consider filtration mappings between graded modules, and that further simplifies the presentation.
 
 Keep in mind that in the decreasing context a  non-zero term $\,\ell_p \in L_p\,$ is the {\it leading term} of an element $\,\ell \in L\,$ if
 \begin{equation} \label{nasgrfive}
 	\ell = \ell_{p} + \cdots + \ell_{ r} + \cdots + \ell_{s} \quad \text{where} \quad \ell_{ j} \in L_{ j}.
 \end{equation}
 We will again write $\,\ell_p\,$ as $\,\overline{\ell}$.

 %\newpage
 %The goal is to develop an algorithm which, for a given non-zero $\,m_p \in M_p,$  \vspfv
 %\begin{itemize}
 %\ibul  determines if   $\,m_p\,$ is the initial term of an element $\,m \in\ker f$,\, and \vspfv
 %\ibul enables one to produce all such $\,m$.  \vspfv
 %\end{itemize}
 For any non-negative integers $\,p\,$ andr $\,r\,$ define \medskip
 \begin{equation} \label{nzdefs}
 	\begin{array}{lllll}
 		Z^p_r & := 
 		& \left\{ \ \ \  
 		\begin{array}{lllllll}
 			\{\,m \in F^{p}M\,:\,fm \in F^{p+r}N\,\} \ \ &\text{if} \ \ r \geq 0, \vspfv \\
 			
 		\end{array} \right.  \quad
 		\vspfv \\
 		M^p_r & := & \Big\{\ \ \ \ \ \big\{\,m_p \in M_p\,:\, m_p = \overline{m}\, \ \ \text{for some}\ \  \ m \in Z^p_r\big\} .\vspfv \\
 	\end{array}
 \end{equation}

  \nd 	Define $\, E_0^{p,0}= {\rm Gr}^pM,  E_0^{p,1}= {\rm Gr}^p N$ and $d_0={\rm Gr}(f)$.  	For $r\geq 1$	define $E_r^{p,0}=M^p_r$ and $E_r^{p,1}=\frac{F^{p}N}{f(Z_{r-1}^{p-r+1}) +  F^{p+1}N},$
 
 \nd  The restriction $\,f\big|_{Z^r_p}\,$ induces a linear mapping
 \begin{equation} \label{thedrmapD}
 	d_r:m_p \in E_r^{p,0} \mapsto [fm] \in E_r^{p+r,1} 
 \end{equation}
 where $\,m \in Z_r^p\,$ is any element satisfying $\,\overline{m} = m_p$.
 
 There are diagrams for the spectral sequence of a decreasing filtration  analogous to \eqref{grdiag}-\eqref{grdiagz}.
 
 \newpage 
 
 \begin{equation}\mbox{Here is the picture for } (E_1,d_1) \end{equation}
 
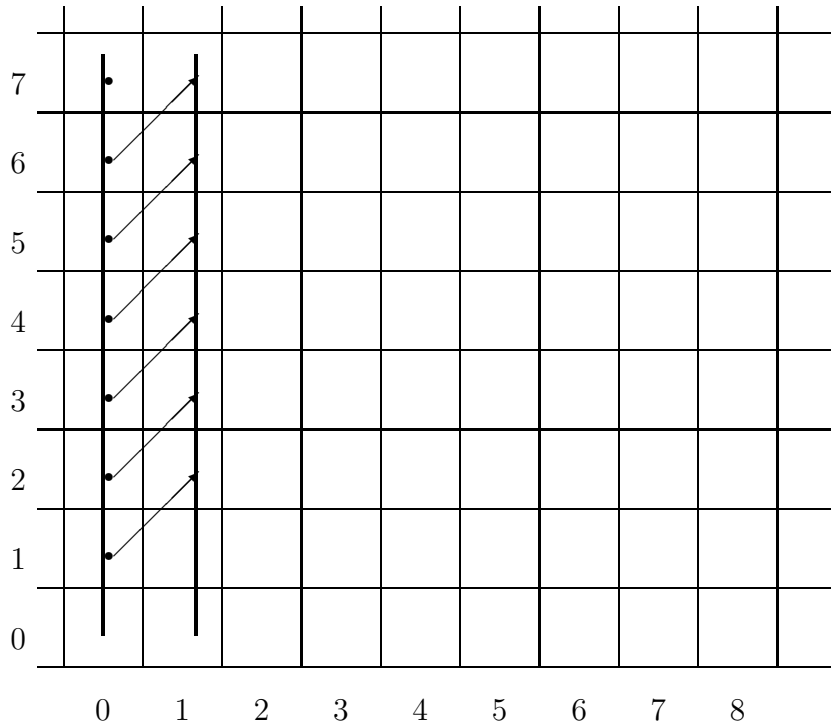
\begin{center}
 	\begin{picture}(420,500)  
 		\hspace*{-.6in}
 		\def\mp{\multiput}
 		\def\elt{\circle*{3}}
 		\mp(89,237)(0,30){6}{\vector(1,1){32}}
 		
 		\mp(60,195)(0,30){9}{\line(1,0){300}}
 		\mp(70,195)(30,0){10}{\line(0,1){250}} \put(142,175){$2$}
 		
 		\put(172,175){$3$} \put(202,175){$4$} \put(232,175){$5$}
 		\put(262,175){$6$} \put(292,175){$7$} \put(322,175){$8$}
 		\put(112,175){$1$} \put(82,175){$0$}
 		
 		\put(50,412){$7$}\put(50,382){$6$}\put(50,352){$5$}
 		\put(50,322){$4$}\put(50,292){$3$}\put(50,262){$2$}
 		\put(50,232){$1$}\put(50,202){$0$}

 		\linethickness{0.35mm}
 		\mp(85,427)(35,0){2}{\line(0,-1){220}}

 		\mp(87,417)(35,0){1}{\elt} \put(120,417)
 		
 		\mp(87,387)(35,0){1}{\elt} \put(120,387)
 		
 		\mp(87,357)(35,0){1}{\elt} 
 		
 		\mp(87,327)(35,0){1}{\elt} 
 		
 		\mp(87,297)(35,0){1}{\elt}

 		\mp(87,267)(35,0){1}{\elt}

 		\mp(87,237)(35,0){1}{\elt}
 		
 		\put(55,155){{\it The } $E_1^{p,q}$ {\it term with differentials from coordinate } $(p,0)$}
 		\put(120,140)	{ {\it to coordinate } $ (p+1,1)$}
 	\end{picture}

 \end{center}

 The following propositions and algorithm are stated for an increasing filtration.  The decreasing case is similar.
 
 \begin{prop} \label{practical}{\rm  For any non-zero $\,m_p \in M_p\,$,   $d_rm_p= [0]$, if and only if there is an element $\,m\pr \in Z^{r+1}_{p}\,$ with leading term $\,m_p$.  \vspfv
 	}
 
 \end{prop} 
 
 To indicate that $\,m_p \in \ker d_r\,$ one says that $m_p$ ``survives'' to  the $r+1$ page of the spectral sequence.   If $\,d_rm_p \neq [0]$ then $\,m_p$ does not represent the leading term of a class in $\text{ker}(f)$.   In other words the $\,A-$modules $\,E_r^{*,1}$  are obstructions to being able to  extend $\, m_p \in M_p$ to $m \in \text{ker}(f)$.  The class $d_rm_p \in E_r^{p-r,1}$ is said to have been ``killed '' by $m_p$.\\

 		\pf

 			 For $\,r  \geq  0\,$ one sees from (\ref{thedrmapD})\, that
 					$d_rm_p = [0]  \Leftrightarrow  \,[fm] = [0] \in \frac{F_{p-r}N}{f(Z^{r-1}_{p-1}) \, + \, F_{p-r-1}N}$ for some $\,m \in Z_p^r $ with leading term $\,m_p$.   This implies there is an $m^{\prime \prime}\in Z_{p-1}^{r-1}$  with $fm = fm^{\prime\prime}$ modulo $F_{p-r-1}N$ .    Then
 					 $\overline{m-m^{\prime\prime}} = m_p$  and $f(m-m^{\prime\prime}) \in  F_{p-r-1}N$,  so we may take $m\pr = m-m^{\prime\prime}$.

 			\qed

 	%\newpage
 	
 	\begin{cor}  \label{keycor}{\rm

 			  $m_p\,$ is the leading term of an element of $ \ \ker f\,$ if and only if $\,m_p$ survives to infinity, i.e.\ if and only if $\,d_rm_p= [0]\,$ for $\,r = 0,1,\dots,p$.
 		}
 	\end{cor}

 	Proposition \ref{practical} and the proof  thereof, together with Corollary \ref{keycor},\, are easily  converted into an algorithm for determining if a given non-zero $\,m_p \in M_p\,$ is the leading term of an element  $\,m \in \ker f$.\\

 \begin{alg} \label{alg}{\rm

 The differentials have the form
 $$d^r: E^r_{p,0} \to E^r_{p-r,1}.$$

 The $E^0_{p,s} \ s=0,1$ terms are given by

 \begin{enumerate}
 	\item $E_{p,0}^0 =F_p M / F_{p-1}M \cong M_p$\\
 	
 	\item $E_{p.1}^0 =F_p N / F_{p-1}N \cong N_p$\\
 	
 	\item $d^0 $ is induced by $f$.
 \end{enumerate}

 \begin{equation*}\tag{i}  \end{equation*} With $ [m] $ denoting the equivalence class in a quotient, an element  $[m_p] \in E_{p,0}^0 =M_p$  survives to  $ E^1_{p,0} $  if $ f(m_p) =n_{p-1}+ \cdots +n_0. $  A class $ [n_p] \in E_{p,1}^0= N_p $ survives to a non-zero class in   $E^1_{p,0}$ if $n_p $  is not in the image of $f|_{M_p}.$

\begin{equation*}\tag{ii}\end{equation*} 
 To compute $d^1([m_p]), $  with $[m_p] \in E^1_{p,0}$ we note the $f(m_p) = n_{p-1}  + \cdots +n_0\in F_{p-1}N$.  Then $d^1[m_p]=[n_{p-1}]$.   If $[n_{p-1}] \neq 0$ in $E^1_{p-1, 1}$ then  $[m_p] $ does not survive to $E^2$ and $m_p$ cannot be the leading term of a class in $\text{ker}(f)$.

\nd  If $[n_{p-1}] = 0$ in $E^1_{p-1,1}$ then $[m_p]$ survives to $E^2_{p,0}$.  To be zero in $E^1_{p-1,1}$ means     $n_{p-1} = f(m_{p-1}) \in f(Z_{p-1}^0)$ modulo $F_{p-2}N$.   The element $[m_p-m_{p-1}]  = [m_p]$ in $E^1_{p,0}$ (i.e. they have the same leading terms) and $[f(m_p-m_{p+1}) ]\in F_{p-2}$ represents $d_2([m_p]) \in E^2_{p-2,1}$.

 \begin{equation*}\tag{iii}\end{equation*} 
 In general $[m_p]$ survives to $E^r_{p,0}$ if $m_p$ can be extended to a class, $m_p + m_{p-1} + \cdots$ such that $f(m_p + m_{p-1} + \cdots) = n_{p-r} + n_{p-r-1} + \cdots$.   The $r-th$ differential,  $d_r([m_p])$ is given by  $ [n_{p-r}]\in E^r_{p-r,1} =\frac{ F_{p-r}N}{ f(Z_{p-1}^{r-1}) + F_{p-r-1}N}$.   The differential  $d_r([m_p]) = 0 \Leftrightarrow  $ if there is such an extension,  $m_p + m_{p-1} + \cdots $ such that  there is $m\pr_{p-1} + m\pr_{p-2} + \cdots $ with $f(m\pr_{p-1} + m\pr_{p-2} + \cdots )=n_{p-r}$ modulo $F_{p-r-1}N$.  The $r-$th differential is given by $[f(m_p +(m_{p-1}-m_{p-1}\pr) + \cdots )]$.   One can check that the differential is independent of choices.

 The primary applications of our spectral sequence  in this paper involve  graded Lie algebras  $\,\cl = \prod_{n \geq 0}\cl_n\,$  over  $A$.
 Elements $\,\ell \in \cl\,$ are written as $\,\ell = \ell_0 + \ell_1 + \ell_2 + \cdots$,\, where $\,\ell_n \in \cl_n$. 
 Elements of $\,\cl_j\,$ are said to be {\it at} {\it level} $\ j$.   The  ``graded'' hypothesis in this context  includes the added assumption that

 \begin{equation} \label{gradcond}
 	[\ell_i,m_j] \in \cl_{i+j} \quad \text{if} \quad \ell_i \in \cl_i \ \ \text{and} \ \ m_j \in \cl_j.
 \end{equation}
 For $\,\ell, m \in \cl\,$ we write the bracket $\,[\ell, m]\,$ as $\,[\ell,m]_0 + [\ell,m]_1 + [\ell,m]_2 + \cdots$,\, where $\,[\ell,m]_j \in \cl_j$.

 For any $\,\ell \in \cl\,$ the  adjoint representation $\,{\rm ad}(\ell):\cl \rightarrow \cl\,$ of $\,\ell$,\, which we refer to informally as the ``ad map of $\,\ell$,''  is defined by 
 \begin{equation} \label{defofad}
 	{\rm ad}(\ell):m \in \cl \mapsto [\ell,m] \in \cl.
 \end{equation}

 \nd It follows from (\ref{gradcond}) that
 \begin{equation} \label{admap}
 	\ell \in \cl_i \quad {\rm ad}(\ell)\big|_{\cl_j}:\cl_j \rightarrow \cl_{i+j}.
 \end{equation}
 
 These $ad$ mappings will be our primary examples of  graded linear maps.   At this level of generality elements of $\,\ker {\rm ad}(m)$,\, which by the skew-symmetry of the bracket include $\,m$,\, are often called {\it centralizers}\, of $\,m$.  We will be interested in detecting elements within their kernels.

 \begin{equation} \label{666}
 	{\rm ad}(\ell_p): m_j \in \mcl_j\mapsto [\ell_p,m_j] \in \mcl_{p+j}.
 \end{equation}

 We apply spectral sequences  to two problems involving a given element
 \begin{equation} \label{thisism}
 	\ell = \ell_0 + \ell_1 + \cdots + \ell_p \in \mcl.
 \end{equation}
 \begin{itemize}
 	\ia Assuming $\,p > 0\,$ give a sufficient condition on $\,{\rm ad}(\ell_p)\,$ for determining all centralizers of $m$:  and 
 	\ib    give a sufficient condition on $\,{\rm ad}(\ell_0)\,$ for determining all formal centralizers, $m$  of $\,\ell$;\,  where a formal centralizer lies in  $\,\prod_{j \geq 0}\mcl_j$ \, rather than  $\,\bigoplus_{j \geq 0}\mcl_j$. 
 \end{itemize}
 
 \

 The conditions we offer and the illustrating examples are elementary.   However, as will be seen in the body of this paper, the hypotheses can be difficult to verify in practice.
 
 We first deal with (b).   Let $M=N=\mcl$ with the decreasing filtration induced by the grading on $\mcl$.     The operator $\text{ad}(\ell_0)$ defines a filtration preserving linear map $M \to N$ which in turn induces a spectral sequence $E_r^{p,q}, \ q=0,1$.  
 
 \begin{prop} \label{viewone} Suppose $\,\ell\,$ is as in\, {\rm(\ref{thisism})}\, and $\,{\rm ad}(\ell_0)\big|_{M_j}:M_j \rightarrow N_j\,$ is surjective for all $\,j > 0$.  Then any $\,m_p\in \ker {\rm (ad}(\ell_0):E_0^{p,0} \to E_0^{p,1} ),$ can be extended to a formal centralizer
 	$$m_p+ m_{p+1} + \cdots $$
 	of $\,\ell$.
 \end{prop}

 \pf  
 At level $\,0\,$ we are given that 
 \begin{equation*}
 	{\rm ad}(\ell_0):\text{Gr}^p M  \to \text{Gr}^p  N
 \end{equation*}
 is surjective.   As a consequence $E_r^{q,1}=0$ for $r>0$ and all higher differentials are $0$.  Any $m_p \in \text{ker} (ad(\ell_0))$ survives to $E_1^{p,0}$.  Since all higher differentials are $0$ such an $m_p$ must extend to a formal centralizer.
 
 \qed \bigskip

 The analogue of Proposition \ref{viewone}\, for\, (a)\,  follows.  We rescale $M$ in order to  have $\text{ad}(\ell_p)$ a filtration preserving map.  To this end   define $M_j=\mcl_{j-p}$ and $N_j=\mcl_j$.  Then $\text{ad}(\ell_p): M_j \to N_j$.  There is a spectral sequence induced by the associated decreasing filtration.
 
 \begin{prop} \label{viewtwo}{\rm  Suppose $\,\ell\,$ is as in {\rm (\ref{thisism})}\, and $\,{\rm ad}(\ell_p)\big|_{M_j}:M_j \rightarrow N_{j}\,$ is surjective for all $\,0 \leq j < q$.   Then any $\,m_q \in \ker {\rm ad}(\ell_p): E^{q,0}_0\to  E^{q,1}_0$ is the leading term of a centralizer of $\,\ell$.  
 }
 \end{prop}
 \qed

}
 
 \end{alg}

 We provide a few elementary examples to show  how the spectral sequence  describe in this appendix works.   For our first example we  consider the polynomial algebra on $n$ variables $P\colon=\mbk[\mathbf{x}]$, $\mathbf{x}=(x_1, \cdots, x_n)$ over a field $\mbk$.   One wants a criterion for a polynomial to belong to an ideal $\mathcal{I}$. A  Gr\"{o}bner basis is a  collection of polynomials  $G=\{g_1,\cdots, g_s\}$ which solves this problem (see  \cite{BoCh} and \cite{clo} ).   Properties familiar   to people who work with Gr\"{o}bner basis may be translated into the language of a spectral sequence associated to a linear map.  Specifically a Gr\"{o}bner basis depends on an ordering on the set of monomials in $\mathbf{x}$  which induces a increasing filtrations on $P$ and $P/\mathcal{I}$.   There is the spectral sequence associated to  the canonical surjection
 $$f:P \longrightarrow P/\mathcal{I}.$$

 The first page of the spectral sequence is   $\{E^0_{p,q}, d^0\}$, where $E^0_{p,q}$ are the associated graded groups of $P$ (for $q=0$)  and $P/\mathcal{I}$ (for $q=1$).  The differential, $d^0$ is  the induced map
 $$d^0=\text{Gr}_p(f):	\text{Gr}_p P  \longrightarrow 	\text{Gr}_p P/\mathcal{I}.$$

It is quite easy to compute  $\text{Gr}_p P $.   The associated graded module induced by the monomial ordering has as most one element in each grading.    Computing  $\text{Gr}_p P/\mathcal{I}$ is more difficult.  The $\mbk$ module  $\text{Gr}_p P/\mathcal{I}$ is  not a new object   to people who work with Gr\"{o}bner basis.   There is the ideal, $\langle \text{Lt} (\mathcal{I}) \rangle$ of leading terms of $\mathcal{I}$ defined as
$$\text{Lt}(\mathcal{I})= \{ cx^a | \text{there exist } f \in I \text{with } \overline{f}=cx^a\},$$
\nd where $\overline{f}$ is the leading term in $f$.  Then

$$\text{Gr}_p P/\mathcal{I} = P/\langle Lt(\mathcal{I} )\rangle.$$

A Gr\"{o}bner basis,   $G=\{g_1,\cdots, g_s\}$  provides a set of polynomial generators of $\mathcal{I}$  that allow a computation of   $ \text{Gr}_p P/\mathcal{I} $.  Specifically $G$ has the property that $\langle Lt(\mathcal{I})\rangle = \langle \{\overline{g}_1, \cdots , \overline{g}_s\}\rangle $.  In particular a $\mbk$ linear basis for 
$ \text{Gr}_p P/\mathcal{I} $ is given by
$$\{u | u  \text{ a monomial in } P \ni u \neq a \overline{g}, g\in G \}.$$

 The connection  between our spectral sequence and Gr\"{o}bner bases is now evident.

 	\begin{itemize} 
 		\ia $m \in \mci$ then $\overline{m}\,$ is divisible by $\,\overline{g}\,$ for some $\,g \in G$.  
 		
 		\nd This is simply the statement that $\overline{m}$ must be in $\text{ker}d^0$.\\
 		
 		\ib Any $\overline{m} \in \text{ker}d^0$ is the leading term of some polynomial in $\mathcal{I}$.
 		
 		\nd  This is obvious since $\overline{m}=hg_i$ for some $i$.  It is also  obviously  a consequence of the spectral sequence since the surjectivity of $d^0$ implies $E^r_{q,1}=0$ for $r>0$.

 		\ic  Given $\overline{m} \in \text{ker}d^0$ the algorithm \eqref{alg} determines a polynomial $m = \sum^n_{i=1} h_ig_i$,\, where $\,h_i \in P\,$ for each $\,i\,$  and $\overline{m}=max_{1\leq i\leq n}\{h_ig_i\}$.
 		
 	\end{itemize}

  The next examples illustrate  algorithm \eqref{alg}.   We first consider  linear ordinary differential  equations.

 \begin{exams} \label{ssexamsone} {\rm  We view solving the differential equation 
 		\begin{equation*} \tag{i}
 			p_n(x)y^{(n)} + p_{n-1}(x)y^{(n-1)} + \cdots + p_0(x)y = 0
 		\end{equation*}
 		with polynomial coefficients    as finding elements of the kernel
 		of the real linear operator $\,f:\nr[x] \rightarrow \nr[x]\,$ given by
 		\begin{equation*} \tag{ii}
 			\ts f:q(x) \mapsto \sum^n_{j=0}p_j(x)q^{(j)}(x).
 		\end{equation*}
 		We assume the increasing  filtrations as in\, \eqref{asgrinc}   but grade the domain and codomain in different ways so as to render $\,f\,$ filtration preserving.   Specifically we define a shift parameter associated to a differential equation.

 		\begin{defin}  {\rm 
 				The {\it shift},  $sh(P)$ of the differential equation 
 				
 				$$P=	p_n(x)y^{(n)} + p_{n-1}(x)y^{(n-1)} + \cdots + p_0(x)y $$
 				
 				\nd is  $\text{max } \{ \text{deg} (p_j)-j\}$.
 				
 			}
 		\end{defin}

 		A term $\,kx^n\,$ in the domain is assigned grade $\,n+sh(P)$,\, whereas in the codomain the same term is graded $\,n$.  When $\,kx^n\,$ is the leading term of some given element, $\,k\,$ would be called the {\it leading coefficient}\, of that element and $\,x^n\,$ would be called the {\it leading monomial}.
 		\begin{itemize}
 			
 			%\newpage
 			
 			\ia  We illustrate algorithm \eqref{alg} by solving the differential equation\\
 			\begin{equation*}  \hspace*{.5in}
 				\ts \big(x^6 - \frac{5}{12}x^3 + \frac{9}{8}x\big)y^{\prime \prime \prime} -
 				\big(6x^5 - \frac{5}{4}x^2 + \frac{9}{8}\big) y\dpr + 15x^4y\pr-15x^3y = 0.
 			\end{equation*}\\ 
 			The shift is 3 and  finding all solutions is equivalent to finding a vector space basis for the kernel of the real linear operator   $\,f:\nr[x] \rightarrow \nr[x]\,$ defined by
 			\begin{equation*}
 				\hspace*{.56in} \ts q \in \nr[x] \mapsto \big(x^6 - \frac{5}{12}x^3 + \frac{9}{8}x\big)q^{\prime \prime \prime} -
 				\big(6x^5 - \frac{5}{4}x^2 + \frac{9}{8}\big) q\dpr + 15x^4q\pr-15x^3q \in \nr[x].
 			\end{equation*}
 			
 			\nd	To apply the algorithm to construct a basis of solutions we need little more than the table \medskip
 			\begin{equation*}\hspace*{.2in}
 				\begin{array}{cccccl}
 					line \ \ & monomial \ \ & grade \  level & f \ o$\!$f & filtration \ level \\
 					&   &   & monomial  \medskip \\
 					1 & x^0 = 1 &   3 & -15x^3 & 3 &  \\
 					2 & x & 4 & 0  & 0\\
 					3 & x^2 & 5 & 3x^5 + {5}{2}x^2 - \frac{9}{4} & 5 \\
 					4 & x^3 & 6 & 5x^3 & 3 \\
 					5 & x^4 &  7 & -3x^7 + 5x^4 + \frac{27}{2}x^2 & 7 \\
 					6 & x^5 & 8 & 45x^3 & 3 \\
 					7 & x^6 & 9 &  15x^9 - \frac{25}{2}x^6 + \frac{405}{4}x^4 & 9 \\
 					8 & x^7 & 10 & 48x^{10}-35x^7+189x^5 & 10 \\
 					9 & x^8 & 11 & 102x^{11}-70x^8+315x^6 & 11 \\
 					10 & x^9 &  12 & 192x^{12}-120x^9+486x^7 & 12 \\
 					%x^{10} & 13 & 315x^{13}-\frac{375}{2}x^{10}+\frac{2834}{4}x^8 & 13
 				\end{array} \smallskip
 			\end{equation*}

 			Among the listed monomials only$\,x^1$,  $\,x^3\,$ and $\,x^5\,$ drop filtration, i.e.\ satisfy $\,d_0x^n = [0]$,\, and are therefore the only possibilities for leading terms of solutions. This takes care of $\,d_0$.     From line 2  $\,y = x\,$ is a solution  and reduces the calculations of $\,d_1x^n\,$ to $\,n = 3\,$ and $\,n = 5$.  	We begin with $\,x^3$.  From the chart we see that $f(x^3)$ drops $3$ filtrations, so represents an element in $E^3_{6,0}$.

 			We see that $\,d_3x^3 = [0]\in E^3_{3,1}$ if there is $kx^m,\, m<3$ with $f(kx^m)=f(x^3)$ modulo filtration $2$.   From the chart we see that $f(-\frac{1}{3}x^0) = f(x^3)$ and that  $\,m = x^3+\frac{1}{3}$ represents the lift of $[x^3]$ to $E^4_{6,0}$, i.e. a representative of $[x^3]$ with image  in filtration $2$  However  $\,fm = 0$ and a second solution has therefore been found.
 			
 			The algorithm for extending $x^5\,$ to solutions is slightly more interesting.  In this case we can see from line 6 and line 1,\, or from line 6 and line 4,\, that $\,d_1x^5 = [0] \in E^1_{7,1}$.  This leads to two distinct solutions with leading term $\,x^5$,\, i.e.\ $\,x^5-9x^3\,$ and $\,x^5+3$.

 			The algorithm has uncovered four solutions, but only three linearly independent solutions are needed for a basis,  e.g.\  $\,\{x,x^3+\frac{1}{3},\,x^5+3\}$. 
 			
 			\ib For the equation
 			\begin{equation*}
 				y\dpr - y = 0
 			\end{equation*}
 			The shift is $0$ and the corresponding operator
 			\begin{equation*}
 				q \in \nr[x] \mapsto q\dpr - q
 			\end{equation*}
 			has the table 
 			\begin{equation*}\hspace*{-.04in}
 				\begin{array}{cccccl}
 					line \ \ & monomial \ \ & grade \  level & f \ o$\!$f & filtration \ level \\
 					&   &   & monomial  \medskip \\
 					1 & x^0 = 1 &   0 & -1 & 0 \\
 					2 & x & 1 &- x  & 1 \\
 					3 & x^2 & 2 &- x^2+2 & 2 \\
 					4 & x^3 & 3 &- x^3 +6x & 3 \\
 					\vdots & \vdots  & \vdots   &  \vdots &  \\
 					n+1 & x^n & n & -x^n + n(n-1)x^{n-2} & n \\
 					\vdots & \vdots  & \vdots   &  \vdots &  
 				\end{array} \smallskip
 			\end{equation*}
 			One sees that $d_0\,$ is an isomorphism for all $\,n \geq 0$.   Hence Algorithm \ref{alg}\, implies there are no  polynomial solutions.  Of course there are no such solutions since a basis is given by $\,\{e^x,\,e^{-x}\}$.
 		\end{itemize}
 	}
 \end{exams}

\begin{exam}{\rm

  Let $\,\mcl\,$ be the LIe algebra  of all upper-triangular $\,n \times n\,$ matrices with real (or complex) coefficients. This becomes a Lie algebra when the bracket $\,[A,B]\,$ is assumed the usual commutator  $\,AB-BA$,\, and becomes a graded Lie algebra by defining $\,L_j := 0\,$ if $\,j > n,\,$ and otherwise defining $\,L_j\,$ to be all those $\,n \times n\,$ matrices $\,[a_{ik}]\,$ having their only non-zero entries, if any, on the $\,(k-i)^{\rm th}$-superdiagonal, i.e.\ $\,k-i\,$ positions above the diagonal.\\

 		To illustrate\, Theorem \ref{viewone}\,  and algorithm \eqref{alg} choose $\,n = 4\,$ in the previous paragraph and let
 			\begin{equation*}
 				m = \left[ \begin{array}{crcrc}
 					5 & 1 & \ 4 & -5 \vspfv \\
 					0& -6 & \ 11 & 3 \vspfv \\
 					0 & 0 &\ 2 & 7 \vspfv \\
 					0 & 0 & \  0 & 1 \end{array} \right],
 			\end{equation*}
 	
 			in which case \medskip
 	$$
	 m_0 = \left[ \begin{array}{crcrc}
 					5 & 0 & \ 0 & \ 0 \vspfv \\
 					0& -6 &   \ 0 & 0 \vspfv \\
 					0 & 0 &\ 2 & 0 \vspfv \\
 					0 & 0 & \  0 & 1 \end{array} \right]\cdots etc. \ \ 
 			$$

 			It is a simple matter to check that $\,{\rm ad}(m_0)\big|_{\mcl_j}:\mcl_j \rightarrow \mcl_j\,$ is an isomorphism for $\,j = 1,2,3$,\, and that $\,{\rm ad}(m_0)\big|_{\mcl_0}:\mcl_0 \rightarrow \mcl_0\,$ is the zero mapping. From  Proposition \ref{viewone}\, we see that any diagonal matrix  may be extended to a commutator of $m$.  In fact this is true if we replace $m$ with any matrix with the same diagonal entries as $m$.  To illustrate the algorithm \eqref{alg} we  extend  
 			\begin{equation*}
 				m\pr_0 = \left[ \begin{array}{crcrc}
 					9 & 0 & \ 0 & 0 \vspfv \\
 					0& 8 &   \ 0 & 0 \vspfv \\
 					0 & 0 &\ 7 & 0 \vspfv \\
 					0 & 0 & \  0 & 6 \end{array} \right].
 			\end{equation*}
 		to  a commutator of $m$.
 			
 			To extend  $[m^{\prime}_0] \in E_1^{0,0}$  to $E_2$ we have to find $m_1^{\prime}$ with  $\text{ad}(m_0)(m_1^{\prime}) = \text{ad}(m)(m_0^{\prime})$.  Since $\text{ad}(m_0)(m_0^{\prime})=0$ this is equivalent to solving for $m_1^{\prime}$ in the equation
 			\begin{equation*}  \tag{i}
 				[m_0,m\pr_1] + [m_1,m\pr_0] = [0]
 			\end{equation*}
 			where $\,[0]\,$ denotes the $\,4 \times 4\,$ matrix of zeroes.   By writing $\,m\pr_1\,$ as
 			\begin{equation*}
 				m\pr_1 =  \left[ \begin{array}{crcrc}
 					0& x_1 & \ 0 & 0 \vspfv \\
 					0&  0 &   x_2 & 0 \vspfv \\
 					0 & 0 & 0 & x_3  \vspfv \\
 					0 & 0 & \  0 & 0 \end{array} \right]
 			\end{equation*}
 			this reduces to
 			\begin{equation*}
 				[0]  =  \left[ \begin{array}{ccccc}
 					0&11x_1 -1 & \ 0 & 0 \vspfv \\
 					0&  0 &   -8x_2-11& 0 \vspfv \\
 					0 & 0 & 0 & x_3-7  \vspfv \\
 					0 & 0 & \  0 & 0 \end{array} \right],
 			\end{equation*}
 			yielding
 			\begin{equation*}
 				m\pr_1 =  \left[ \begin{array}{crcrc}
 					0& \frac{1}{11} & \ 0 & 0 \vspfv \\
 					0&  0 &  -\frac{11}{8} & 0 \vspfv \\
 					0 & 0 & \ 0 & 7  \vspfv \\
 					0 & 0 & \  0 & 0 \end{array} \right].
 			\end{equation*}
 			
 			In the language of our spectral sequence $m_0\pr$ survives to $E_2$ and is represented by  $m_0\pr + m_1\pr$.   Continuing in this way, the next step in the algorithm is to find $\,m\pr_2\,$ with
 			\begin{equation*} \tag{ii}
 				[m_0,m\pr_2] + [m_1,m\pr_1] + [m_2,m\pr_0] = [0]
 			\end{equation*}

 		  By writing $\,m\pr_2\,$ as an indeterminate filtration 2 matrix as in the previous step we find 
 			
 			\begin{equation*}
 				m\pr_2  = \left[ \begin{array}{crccc}
 					0& 0 & \ \frac{83}{24} & 0 \vspfv \\
 					0&  0 &  0 & \frac{645}{56} \vspfv \\
 					0 & 0 & 0 & 0  \vspfv \\
 					0 & 0 & \  0 & 0 \end{array} \right].
 			\end{equation*}
 		
 		\nd and  $m_o\pr$ survives to $E_3 $ represented by $m_0\pr+m_1\pr+m_2\pr$.\\
 			
 			It remains to solve the equation
 			\begin{equation*}
 				[m_0,m\pr_3] + [m_1,m\pr_2] + [m_2,m\pr_1] + [m_3,m\pr_0] = [0]
 			\end{equation*}
 			for $\,m\pr_3$.  By writing $\,m\pr_3\,$ as
 			\begin{equation*}
 				m\pr_3 =  \left[ \begin{array}{crccc}
 					0& 0 &0 & x_1 \vspfv \\
 					0&  0 &  0 & 0\vspfv \\
 					0 & 0 & 0 & 0  \vspfv \\
 					0 & 0 &   0 & 0 \end{array} \right],
 			\end{equation*}
 			and following the same procedure used with the first three equations
 			one finds that
 			\begin{equation*}
 				m\pr_3 =  \left[ \begin{array}{crccc}
 					0& 0 &0 & -\frac{13,877}{1,848} \vspfv \\
 					0&  0 &  0 & 0\vspfv \\
 					0 & 0 & 0 & 0  \vspfv \\
 					0 & 0 &   0 & 0 \end{array} \right].
 			\end{equation*}
 			
 			Once we have found a representative of the lift of $m_0\pr$ to $E_4$ we have found a centralizer of $\,m\,$ with initial term $\,m\pr_0\,$:
 			\begin{equation*}
 				m\pr = m\pr_0 + m\pr_1 + m\pr_2 + m\pr_3 = \left[ \begin{array}{crccc}
 					9& \frac{1}{11} & \frac{83}{24}  & -\frac{13,877}{1,848} \vspfv \\
 					0&  8 &  -\frac{11}{8} & \frac{645}{56}\vspfv \\
 					0 & 0 & 7 & 7  \vspfv \\
 					0 & 0 &   0 & 6 \end{array} \right].
 			\end{equation*}
 			In fact one can see from the algebraic steps that this must be the unique centralizer of $\,m\,$ with initial term $\,m\pr_0$.
 	
}
 
\end{exam}

 \section{Hamiltonian Vector Fields}\label{sec:hamiltonian}

 In this appendix we gather information about Hamiltonian vector fields.  The collection $\ch$ of formal power series in the variables $x=(x_1, \cdots,x_{2n})$, assumed to begin with quadratic terms is a real Lie algebra with bracket given by the Poisson bracket.

 \begin{equation} \label{def:poisson}
 \{H,K\} := \langle J H_x, K_x \rangle ,\end{equation}

\nd where $J$ denotes the  canonical  matrix $\left[ \begin{matrix}
0 & I \\
-I & 0 
\end{matrix}\right]$ ($I$ being the $n \times n$ identity matrix), $H_x$ and $K_x$ are the respective formal gradients of $H$ and $K$, regarded as column vectors and $\langle , \rangle$ is the column-vector adaptation of the usual inner product.

The Lie algebra $\ch$ becomes a graded Lie algebra if for each $i \geq 0$ we take $\ch_i \subset \ch$ to be the collection of homogeneous polynomials of degree $i+2$.

When viewed as a vector field $JH_x$ in  \eqref{def:poisson} is generally expressed as $X_H$.  A centralizer $K \in \ch$ of $H$ is known classically as a {\it first integral } of $X_H$.

With this notation Hamilton's equation for a {\it Hamiltonian (function)} $H:\nr^{2n} \to \nr$ becomes

\begin{equation}
\label{def:hameq}
\dot x =J H_x(x) \mbox{ where } \dot x = \frac{dx}{dt}.\end{equation}

The connection between Hamilton's equation and \eqref{def:poisson} is given by:

\begin{equation}\label{chainrule}
	\frac{d}{dt}K(x(t)) = \frac{dK}{dx}\frac{dx}{dt} =\langle K_x, JH_x \rangle  = \{H,K\}.\end{equation}

\nd	From which we see that $K$ is a centralizer of \eqref{def:poisson} if and only if $K$ is constant along solutions, which in turn is the case if and only if $K\in ker(ad(H))$.

When dealing with Hamiltonian systems it is useful to keep the following properties of the canonical matrix in mind:

\begin{equation}
\label{proprties:J}
J^2 = - I \quad J^{\tau}=-J =J^{-1}.
\end{equation}

\begin{equation}\label{propertiesH}\textbf{Properties of Hamiltonians}\end{equation}
	
	\begin{enumerate}[label=(\alph*)]

		\item (Conservation of Energy) $H$ is an integral of $X_H$.

		This follows because
		
		$$\{H,H\} = \langle H_x, J H_x\rangle =0$$

		\item (Functional Dependence)  Suppose $K:\nr^{2n} \to \nr$ is a rational function, i.e. assume $K\in \nr(x_1, x_2, \cdots, x_{2n})$ and suppose there is a rational function $g(s,t) \in \nr(x,t)$ such that

		\begin{enumerate}[label=(\roman*)]
			\item  $$\dfrac{\partial g}{\partial s} \neq 0,$$

			and

			\item  $$g(K,H)=0$$
		\end{enumerate}

		Then $K$ is an integral of $X_H$.

		By (ii) we have $\dfrac{\partial g}{\partial s} K_x  + \dfrac{\partial g}{\partial t} H_x = 0$.  Taking the inner product with $JH_x$ and using (a)  above then gives

		$$0=\dfrac{\partial g}{\partial s} \langle K_x, JH_x\rangle   + \dfrac{\partial g}{\partial t}\langle H_x, JH_x\rangle = \dfrac{\partial g}{\partial s} \langle K_x, JH_x\rangle.$$

\nd		whereupon from (i) we conclude that $\langle K_x, JH_x\rangle=0$.

\item For any non-zero integer, $m$, the function $K = H^m$ is a first integral of $X_H$:  take $g(s,t)=s^m - t$ in (b).

\item \label{bd} If $K$ is a first integral of $X_H$ then for every positive integer $m$

\nd   $$\{H,K^m\} := m K^{n-1}\langle J H_x, K_x \rangle =0.$$.

It follows that $K^m$ is a first integral of $X_H$ as well.

\item Suppose $H:\nr^{2n} \to \nr$ is a Hamiltonian function of the form
$H = \frac{1}{2} \langle x, Ax\rangle$, where $A$ is symmetric and invertible.  Define $B:=-JA^{-1}J$ and $K:\nr^{2n} \to \nr$ by $K=\frac{1}{2} \langle x, Bx \rangle$.  The $K$ is a first integral of $X_H$.

To see this recall that for $Q(x) = \langle x, Tx \rangle$ where $T:\nr^n \to \nr^n$ is a linear operator with usual basis  matrix,  $M$ we have

$$Q_x=\dfrac{dQ}{dx} = (M + M^{\tau})x.$$

In particular since  $A$ and $B$ are both symmetric \newline $H_x = Ax \quad K_x=Bx$ and we have

$$
\begin{array}{llll}
	\langle K_x,JH_x\rangle &=& \langle Bx,JAx\rangle \\
	&=&  \langle x,BJAx\rangle \\
	&=&  \langle x,(-JA^{-1}J)JAx\rangle\\
	&=&  \langle x,JA^{-1}Ax\rangle\\
	&=& \langle x,Jx\rangle\\
	&=&0

	\end{array}$$

\begin{equation}\label{PBrequivalent}
\{H,K\} = - \dfrac{dH}{dx} J ( \dfrac{dK}{dx})^{\tau}, \quad\mbox{ where }  \tau \mbox{ is  the matrix } \mbox{ transpose. }
\end{equation}

This is immediate from:

$$
\begin{array}{lllll}
	\{H,K\} &=& \{JH_x, K_x\} \\
	&=& \langle H_x, - JK_x \rangle \quad \mbox{ ( because } J^{\tau}=-J)\\
	&=&-\langle (\frac{dH}{dx})^{\tau}, J (\frac{dK}{dx})^{\tau} \\
	&=&  - \dfrac{dH}{dx} J ( \dfrac{dK}{dx})^{\tau}.
	\end{array}$$
		
	\end{enumerate}

For each $H \in \ch$ the mapping $ad(H);\ch \to \ch$ is a real linear operator.  In addition we have

\begin{prop}\label{prop:derivation}{\rm 
   The operator $ad(H)$ is a derivation.
}
\end{prop}

We remark that an associative algebra with a  Lie bracket such that   $ad(H)$ is a derivation is called a Poisson algebra.\\

\pf for any solution, $x=x(t)$ of \eqref{def:hameq} one has

$$\begin{array}{llll}
ad(H)(K_1K_2) &=& \{K_1K_2,H\} \\
&=& \frac{d}{dt} ( (K_1K_2)(x(t))) \\
&=& \frac{d}{dt}(K_1(x(t)) K_2(x(t)) \\
&=& K_1(x(t)) \frac{d}{dt}(K_2(x(t)) + \frac{d}{dt}(K_1(x(t)) K_2((x(t))\\
&=& K_1 \cdot ad(H)(K_2) + (ad(H)(K_1)\cdot K_2
\end{array}$$ \qed

As a corollary we have that the integrals of any $H \in \ch$ form  a  subalgebra of $\ch$.

\begin{rmk}{\rm 
\label{eqofmotion}
The Hamiltonian system  \eqref{def:hameq}  is often expressed as

\begin{equation}
   \begin{array}{cccc}
\dot{q}_i &= &\dfrac{\partial H}{\partial p_i} \\
\dot{p}_i &=& -\dfrac{\partial H}{\partial q_i}
\end{array}\tag{i}
\end{equation}

Or more succinctly as

\

\begin{equation}  \label{hamcanonical}
\begin{array}{cccc}
\dot{q} &= & H_p\\
\dot{p} &=& -H_q
\end{array}\tag{ii}
\end{equation}

\nd where 

\begin{equation}\label{canonical variables}
x=(x_1,x_2,\cdots, x_{2n}) =: (q_1,q_2,\cdots,q_n,p_1,p_,\cdots, p_n).\tag{iii}
\end{equation}

The $q_i$ and $p_i$ are called {\it canonical variables}.  The $q_i$ are the {\it position variables} and the $p_i$ are the {\it momentum variables}.

An important special case occurs when $H=H(q,p) = \frac{1}{2} |p|^2 +V(q)$.  For example the quadratic term of the Henon-Heiles Hamiltonian \eqref{def:henon-heiles}.  In this context equations \eqref{hamcanonical} become

\begin{equation}  \label{Newton}
\begin{array}{cccc}
\dot{q} &= & p\\
\dot{p} &=& -V_q
\end{array}
\end{equation}
}
\end{rmk}

\nd which are obviously equivalent to Newton's equations
$$   \ddot{q}= -V_q$$

\nd for a particle  of unit mass moving in a conservative force field with potential $V$.

 \

\bibliographystyle{amsalpha}

\end{document}